\documentclass{article}%
\usepackage{amssymb}
\usepackage{amsmath}
\usepackage{graphicx}
\usepackage{amsfonts}%
\providecommand{\U}[1]{\protect\rule{.1in}{.1in}}
\begin{document}

\title{A Theory of Adjoint Functors ---with some Thoughts about their Philosophical Significance}
\author{David Ellerman\thanks{This paper is dedicated to the memory of Gian-Carlo
Rota---mathematician, philosopher, mentor, and friend.}}
\date{}
\maketitle

\begin{abstract}
The question \textquotedblleft What is category theory\textquotedblright\ is
approached by focusing on universal mapping properties and adjoint functors.
Category theory organizes mathematics using morphisms that transmit structure
and determination. Structures of mathematical interest are usually
characterized by some universal mapping property so the general thesis is that
category theory is about determination through universals. In recent decades,
the notion of adjoint functors has moved to center-stage as category theory's
primary tool to characterize what is important and universal in mathematics.
Hence our focus here is to present a theory of adjoint functors, a theory
which shows that all adjunctions arise from the birepresentations of
\ \textquotedblleft chimeras\textquotedblright\ or \textquotedblleft
heteromorphisms\textquotedblright\ between the objects of different
categories. Since representations provide universal mapping properties, this
theory places adjoints within the framework of determination through
universals. The conclusion considers some unreasonably effective analogies
between these mathematical concepts and some central philosophical themes.
[Forthcoming in: \textit{What is Category Theory?} Giandomenico Sica ed.,
Milan: Polimetrica.]

\end{abstract}
\tableofcontents

\section{Introduction: What is Category Theory?}

How might a question like \textquotedblleft What is category
theory?\textquotedblright\ be approached? Technically, the answer is
well-known so that cannot be the point of the question.{\footnote{Some
familiarity with basic category theory is assumed. Whenever possible, I will
follow MacLane \cite{mac:cwm} on notation and terminology. Proofs will be
avoided in favor of citations to the literature unless \textquotedblleft the
pudding is in the proof.\textquotedblright}} The sense of the question is more
about the philosophical or foundational importance of category theory (CT).

One proposed answer might be that CT provides the language in which to
formulate topos theory which, in turn, provides a massive
topologically-flavored generalization of the set theoretic foundations of
mathematics. Although topos theory has been of great importance to release set
theory's `death grip' on foundations (and although my 1971 dissertation
$\cite{ell:sh}$ was on generalizing ultraproducts to sheaves), I do not
believe that the foundational importance of CT is exhausted by providing
generalizations of set theory.

To understand category theory in its own way, we must return to the basic idea
that CT sees the world of mathematics not just in terms of objects but in
terms of morphisms. Morphisms express the transmission of structure or,
broadly, the transmission of determination horizontally\ between objects. The
morphism-view of the mathematical world gives quite a different perspective
than set theory.

Both theories involve universals. But there are two types of universals, the
non-self-participating, vertical, or \textquotedblleft
abstract\textquotedblright\ universals of set theory and the
self-participating, horizontal, or \textquotedblleft
concrete\textquotedblright\ universals which are given by the universal
mapping properties (UMPs) of category theory (see Ellerman $\cite{ell:cu}%
,\cite{ell:it}$). The abstract or non-self-participating universals of set
theory collect together instances of a property but involve no machinery about
the determination of the instances having the property. In contrast, a
self-participating universal has the property itself and every other instance
of the property has it by participating in (e.g., uniquely factoring through)
the universal. The determination that the other instances have the property
`flows through' the universal to the instances.

In brief, our proposed answer to the question of \textquotedblleft What is
CT?\textquotedblright\ is that category theory is about determination
represented by morphisms and the central structure is determination through
the universals expressed by the UMPs.

Universals seem to always occur as part of an adjunction. Hence this research
programme leads from the \textquotedblleft What is CT?\textquotedblright%
\ question to a focus on adjoint functors.\ It is now widely recognized that
adjoint functors characterize the structures that have importance and
universality in mathematics. Our purpose here is to give a theory of
\textquotedblleft what adjoint functors are all about\textquotedblright\ that
will sustain and deepen the thesis that category theory is about determination
through universals.

Others have been lead to the focus on adjoints by different routes. Saunders
MacLane and Samuel Eilenberg famously said that categories were defined in
order to define functors, and functors were defined in order to define natural
transformations. Their original paper $\cite{eilmac:gte}$ was named not
\textquotedblleft General Theory of Categories\textquotedblright\ but
\textit{General Theory of Natural Equivalences}. Adjoints were (surprisingly)
only defined later $\cite{kan:af}$ but the realization of their foundational
importance has steadily increased over time $\cite{law:adj,lam:her}$. Now it
would perhaps be not too much of an exaggeration to see categories, functors,
and natural transformations as the prelude to defining adjoint functors. As
Steven Awodey put it in his (forthcoming) text:

\begin{quote}
{\footnotesize The notion of adjoint functor applies everything that we've
learned up to now to unify and subsume all the different universal mapping
properties that we have encountered, from free groups to limits to
exponentials. But more importantly, it also captures an important mathematical
phenomenon that is invisible without the lens of category theory. Indeed, I
will make the admittedly provocative claim that adjointness is a concept of
fundamental logical and mathematical importance that is not captured elsewhere
in mathematics. }$\cite{aw:cat}$
\end{quote}

Other category theorists have given similar testimonials.

\begin{quote}
{\footnotesize To some, including this writer, adjunction is the most
important concept in category theory. \cite[p. 6]{wood:ord}}

{\footnotesize Nowadays, every user of category theory agrees that
[adjunction] is the concept which justifies the fundamental position of the
subject in mathematics. \cite[p. 367]{tay:pfm}}
\end{quote}

\noindent Hence a theory of adjoint functors should help to elucidate the
foundational importance of category theory.

\section{Overview of the Theory of Adjoints}

It might be helpful to begin with a brief outline of the argument. The basic
building blocks of category theory are categories where a category consists of
objects and (homo-)morphisms between the objects \textit{within} the category,
functors as homomorphisms between categories, and natural transformations as
morphisms between functors. But there is another closely related type of
entity that is routinely used in mathematical practice and is not `officially'
recognized in category theory, namely morphisms directly between objects in
different categories such as the insertion of the generators $x\Rightarrow Fx$
from a set into the free group generated by the set. These cross-category
object morphism will be indicated by double arrows $\Rightarrow$ and will be
called \textit{chimera morphisms} (since their tail is in one category, e.g.,
a set, and their head is in another category, e.g., a group) or
\textit{heteromorphisms} in contrast with homomorphisms.

Since the heteromorphisms do not reside within a category, the usual
categorical machinery does not define how they might compose. But that is not
necessary. Chimera do not need to `mate' with other chimera to form a
`species' or category; they only need to mate with the intra-category
morphisms on each side to form other chimera. The appropriate mathematical
machinery to describe that is the generalization of a group acting on a set to
a generalized monoid or category acting on a set (where each element of the
set has a \textquotedblleft domain\textquotedblright\ and a \textquotedblleft
codomain\textquotedblright\ to determine when composition is defined). In this
case, it is two categories acting on a set, one on the left and one on the
right. Given a chimera morphism $c:x\Rightarrow a$ from an object in a
category $\mathbf{X}$ to an object in a category $\mathbf{A}$ and morphisms
$h:x^{\prime}\rightarrow x$ in $\mathbf{X}$ and $k:a\rightarrow a^{\prime}$ in
$\mathbf{A}$, the composition $ch:x^{\prime}\rightarrow x\Rightarrow a$ is
another chimera $x^{\prime}\Rightarrow a$ and the composition $kc:x\Rightarrow
a\rightarrow a^{\prime}$ is another chimera $x\Rightarrow a^{\prime}$ with the
usual identity, composition, and associativity properties. Such an action of
two categories acting on a set on the left and on the right is exactly
described by a het-bifunctor $Het:\mathbf{X}^{op}\times\mathbf{A}%
\rightarrow\mathbf{Set}$ where $Het(x,a)=\{x\Rightarrow a\}$ and where
$\mathbf{Set}$ is the category of sets and set functions. Thus the natural
machinery to treat object-to-object chimera morphisms between categories are
het-bifunctors $Het:\mathbf{X}^{op}\times\mathbf{A}\rightarrow\mathbf{Set}$
that generalize the hom-bifunctors $Hom:\mathbf{X}^{op}\times\mathbf{X}%
\rightarrow\mathbf{Set}$ used to treat object-to-object morphisms within a category.

How might the categorical properties of the het-bifunctors be expressed
without overtly recognizing chimera? Represent the het-bifunctors using
hom-bifunctors on the left and on the right! Any bifunctor $\mathcal{D}%
:\mathbf{X}^{op}\times\mathbf{A}\rightarrow\mathbf{Set}$ is
\textit{represented on the left}\footnote{This terminology \textquotedblleft
represented on the left\textquotedblright\ or \textquotedblleft on the
right\textquotedblright\ is used to agree with the terminology for left and
right adjoints.} if for each $x$ there is an object $Fx$ in $\mathbf{A}$ and
an isomorphism $Hom_{\mathbf{A}}(Fx,a)\cong\mathcal{D}(x,a)$ natural in $a$.
It is a standard result that the assignment $x\mapsto Fx$ extends to a functor
$F$ and that the isomorphism is also natural in $x$. Similarly, $\mathcal{D}$
is \textit{represented on the right} if for each $a$ there is an object $Ga$
in $\mathbf{X}$ and an isomorphism $\mathcal{D}(x,a)\cong Hom_{\mathbf{X}%
}(x,Ga)$ natural in $x$. And similarly, the assignment $a\mapsto Ga$ extends
to a functor $G$ and that the isomorphism is also natural in $a$.

If a het-bifunctor $Het:\mathbf{X}^{op}\times\mathbf{A}\rightarrow
\mathbf{Set}$ is represented on both the left and the right, then we have two
functors $F:\mathbf{X}\rightarrow\mathbf{A}$ and $G:\mathbf{A}\rightarrow
\mathbf{X}$ and the isomorphisms natural in $x$ and in $a$:
\[
Hom_{\mathbf{A}}(Fx,a)\cong Het(x,a)\cong Hom_{\mathbf{X}}(x,Ga).
\]

\noindent It only remains to drop out the middle term $Het(x,a)$ to arrive at
the wonderful \textit{pas de deux} of the `official' definition of a pair of
adjoint functors---without any mention of heteromorphisms. That, in short, is
our theory of \textquotedblleft what adjoint functors are \textit{really}
about.\textquotedblright\ Adjoint functors are of foundational relevance
because of their ubiquity in picking out important structures in ordinary
mathematics. For such concretely occurring adjoints, the heteromorphisms can
be easily recovered (e.g., in all of our examples). But when an adjunction is
abstractly defined---as always \textit{sans} middle term---then where are the chimeras?

Hence to round out the theory, we give an \textquotedblleft adjunction
representation theorem\textquotedblright\ which shows how, given any
adjunction $F:\mathbf{X}\rightleftarrows\mathbf{A}$ $:G$, heteromorphisms can
be defined between (isomorphic copies of) the categories $\mathbf{X}$ and
$\mathbf{A}$ so that (isomorphic copies of) the adjoints arise from the
representations on the left and right of the het-bifunctor. The category
$\mathbf{X}$ is embedded in the product category $\mathbf{X}\times\mathbf{A}$
by the assignment $x\mapsto(x,Fx)$ to obtain the isomorphic copy
$\widehat{\mathbf{X}}$, and $\mathbf{A}$ is embedded in the product category
by $a\mapsto(Ga,a)$ to yield the isomorphic copy $\widehat{\mathbf{A}}$. Then
the properties of the adjunction can be nicely expressed by the commutativity
within $\mathbf{X}\times\mathbf{A}$ of \textquotedblleft adjunctive
squares\textquotedblright\ of the form:

\begin{center}
$%
\begin{array}
[c]{ccccc}
& (x,Fx) & \overset{(f,Ff)}{\longrightarrow} & (Ga,FGa) & \\
_{(\eta_{x},1_{Fx})} & \downarrow &  & \downarrow & _{(1_{Ga},\varepsilon
_{a})}\\
& (GFx,Fx) & \overset{(Gg,g)}{\longrightarrow} & (Ga,a) &
\end{array}
$
\end{center}

\noindent where the main diagonal $(f,g)$ in a commutative adjunctive square
pairs together maps that are images of one another in the adjunction
isomorphism $g=f^{\ast}$ and $f=g^{\ast}$ (i.e., are adjoint transposes of one
another). Since the maps on top are always in $\widehat{\mathbf{X}}$ and the
maps on the bottom are in $\widehat{\mathbf{A}}$, the main diagonal pairs of
maps (including the vertical maps)---which are ordinary morphisms in the
product category---have all the categorical properties of chimera morphisms
from objects in $\mathbf{X}\cong\widehat{\mathbf{X}}$ to objects in
$\mathbf{A}\cong\widehat{\mathbf{A}}$. Hence the heteromorphisms are
abstractly defined as the pairs of adjoint transposes,
$Het(x,a)=\{(x,Fx)\overset{(f,f^{\ast})}{\longrightarrow}(Ga,a)\}$, and the
adjunction representation theorem is that (isomorphic copies of) the original
adjoints $F$ and $G$ arise from the representations on the left and right of
this het-bifunctor.

Another bifunctor $\mathcal{Z}(Fx,Ga)$ of chimeras is later defined whose
elements are heteromorphisms in the opposite direction (from $\mathbf{A}$ to
$\mathbf{X}$). They are defined only on the images $Fx$ and $Ga$ of the pair
of adjoint functors. These heteromorphisms would be represented by the
southwest-to-northeast anti-diagonal maps $(Gg,Ff)$ in an adjunctive square.
Moreover, this fourth bifunctor is naturally isomorphic to the other three
bifunctors.%
\[
Hom_{\mathbf{A}}(Fx,a)\cong Het(x,a)\cong\mathcal{Z}(Fx,Ga)\cong
Hom_{\mathbf{X}}(x,Ga).
\]

\noindent The universals of the unit $\eta_{x}$ and counit $\varepsilon_{a} $
are associated in the isomorphism respectively with the identity maps $1_{Fx}$
and $1_{Ga}$. The elements of the other two chimera bifunctors associated with
these identities also have universality properties. The identity $1_{Fx}$ is
associated with $h_{x}\ $in $Het(x,Fx)$ and $h_{x2}$ in $\mathcal{Z}(Fx,GFx)$,
and $1_{Ga}$ is associated with $e_{a}$ in $Het(Ga,a)$ and with $e_{a1}$ in
$\mathcal{Z}(FGa,Ga)$. The two $h$-universals provide an over-and-back
factorization of the unit: $x\overset{h_{x}}{\Longrightarrow}Fx\overset
{h_{x2}}{\Longrightarrow}GFx=x\overset{\eta_{x}}{\longrightarrow}GFx$, and the
two $e$-universals give an over-and-back factorization of the counit:
$FGa\overset{e_{a1}}{\Longrightarrow}Ga\overset{e_{a}}{\Longrightarrow
}a=FGa\overset{\varepsilon_{a}}{\longrightarrow}a$. These chimera universals
provide another factorization of any $f:x\rightarrow Ga$ in additional to the
usual one as well as another factorization of any $g:Fx\rightarrow a$. There
is also another over-and-back factorization of $1_{Fx}$ and of $1_{Ga}$ in
addition to the triangular identities. Moreover, there is a new type of
all-chimera factorization. Given any heteromorphism $x\overset{c}{\Rightarrow
}a$, there is a unique $\mathbf{A}$-to-$\mathbf{X}$ chimera morphism
$Fx\overset{z(c)}{\Longrightarrow}Ga$ that factors $c$ through the chimera
version of the unit, i.e., $h_{x}$, \textit{and} through the chimera version
of the counit, i.e., $e_{a}$, in the over-back-and-over-again or zig-zag factorization:%

\[
x\overset{c}{\Rightarrow}a=x\overset{h_{x}}{\Longrightarrow}Fx\overset
{z(c)}{\Longrightarrow}Ga\overset{e_{a}}{\Longrightarrow}a.
\]

\noindent Roughly speaking, the chimeras show their hybrid vigor by more than
doubling the number of factorizations and identities associated with an adjunction.

What may be new and what isn't new in this theory of adjoints? The theory
contains no strikingly new formal results; the level of the category theory
involved is all quite basic. The heteromorphisms are formally treated using
bifunctors of the form $Het:\mathbf{X}^{op}\times\mathbf{A}\rightarrow
\mathbf{Set}$. Such bifunctors and generalizations replacing $\mathbf{Set}$ by
other categories have been studied by the Australian school under the name of
\textit{profunctors }\cite{kelly:enrich}, by the French school under the name
of \textit{distributors} \cite{ben:distr}, and by William Lawvere under the
name of \textit{bimodules} \cite{law:met}.\footnote{Thanks to John Baez for
these connections with the literature on enriched categories.} However, the
guiding interpretation has been interestingly different. \textquotedblleft
Roughly speaking, a distributor is to a functor what a relation is to a
mapping\textquotedblright\ \cite[p. 308]{bor:hca1} (and hence the name
\textquotedblleft profunctor\textquotedblright\ in the Australian school). For
instance, if $\mathbf{Set}$ was replaced by $\mathbf{2}$, then the bifunctor
would just be the characteristic function of a relation from $\mathbf{X}$ to
$\mathbf{A}$. Hence in the general context of enriched category theory, a
\textquotedblleft bimodule\textquotedblright\ $Y^{op}\otimes X\overset
{\varphi}{\longrightarrow\mathcal{V}}$ would be interpreted as a
\textquotedblleft$\mathcal{V}$-valued relation\textquotedblright\ and an
element of $\varphi(y,x)$ would be interpreted as the \textquotedblleft
truth-value of the $\varphi$-relatedness of $y$ to $x$\textquotedblright%
\ \cite[p. 158 (p. 28 of reprint)]{law:met}. The subsequent development of
profunctors-distributors-bimodules has been along the lines suggested by that
guiding interpretation.

In the approach taken here, the elements $x\Rightarrow a$ in $Het(x,a)$ are
interpreted as heteromorphisms from an object $x$ in $\mathbf{X}$ to an object
$a$ in $\mathbf{A}$ on par with the morphisms within $\mathbf{X}$ or
$\mathbf{A}$, not as an element in a `relational' generalization of a functor
from $\mathbf{X}$ to $\mathbf{A}$. Such chimeras exist in the wild\ (i.e., in
mathematical practice) but are not in the `official' ontological zoo of
category theory that sees object-to-object morphisms as only existing
\textit{within} a category. The principal novelty here (to my knowledge) is
the use of the chimera morphism interpretation of these bifunctors to carry
out a whole program of interpretation for adjunctions, i.e., a \textit{theory}%
\ of adjoint functors. In the concrete examples, chimera morphisms have to be
\textquotedblleft found\textquotedblright\ as is done in the broad classes of
examples treated here. However, in general, the adjunction representation
theorem uses a very simple construction to show how `abstract' heteromorphisms
can always be found so that any adjunction arises (up to isomorphism) out of
the representations on the left and right of the het-bifunctor of such heteromorphisms.

These conceptual structures suggest various applications discussed in the
conclusions.\footnote{See Lambek \cite{lam:her} for rather different
philosophical speculations about adjoint functors.} Following this overview,
we can now turn to a more leisurely development of the concepts. Universal
mapping properties and representations are two ways in which \textquotedblleft
universals\textquotedblright\ appear in category theory so we next turn to the
contrasting universals of set theory and category theory.

\section{Universals}

The general notion of a universal for a property is ancient. In Plato's Theory
of Ideas or Forms ($\epsilon\iota\delta\eta$), a property $F$ has an entity
associated with it, the \textit{universal} $u_{F}$, which uniquely represents
the property. An object $x$ has the property $F$, i.e., $F(x)$, if and only if
(iff) the object $x$ \textit{participates} in the universal $u_{F}$. Let $\mu$
(from $\mu\epsilon\theta\epsilon\xi\iota\sigma$ or methexis) represent the
participation relation so \textquotedblleft$x\,\mu\,u_{F}$\textquotedblright%
\ reads as \textquotedblleft$x$ participates in $u_{F}$\textquotedblright.
Given a relation $\mu$, an entity $u_{F}$ is said to be a \textit{universal
for the property} $F$ (with respect to $\mu$) if it satisfies the following
\textit{universality condition}: for any $x$, $x\,\mu\,u_{F}$ if and only if
$F(x)$.

A universal representing a property should be in some sense unique. Hence
there should be an equivalence relation ($\cong$) so that universals satisfy a
\textit{uniqueness condition}: if $u_{F}$ and $u_{F}^{\prime}$ are universals
for the same $F$, then $u_{F}$ $\cong$ $u_{F}^{\prime}$. These might be taken
as the bare essentials for any notion of a universal for a property.

Set theory defines \textquotedblleft participation\textquotedblright\ as
membership represented by the $\in$ taken from $\epsilon\iota\delta\eta$. The
universal for a property is set of objects with that property ${\{x|F(x)\}}$,
the \textquotedblleft extension\textquotedblright\ of the property, so the
universality condition becomes the comprehension scheme:

\begin{center}
$x$ $\in{\{x|F(x)\}}$ iff $F(x)$.
\end{center}

The equivalence of universals for the same property takes the strong form of
identity of sets if they have the same members (extensionality). The
universals of set theory collect together objects or entities with the
property in question to form a new entity; there is no machinery or structure
for the universal to `determine' that the instances have the property. But
`if' they have the property, then they are included in the set-universal for
that property.

The Greeks not only had the notion of a universal for a property; they also
developed the notion of hubris. Frege's hubris was to try to have a general
theory of universals that could be either self-participating or
not-self-participating. This led to the paradoxes such as Russell's paradox of
the universal $R$ for all and only the universals that are not
self-participating. If $R$ does not participate in itself, then it must
participate in itself. And if $R$ does participate in itself then it must be
non-self-participating.{\footnote{Note that Russell's paradox was formulated
using the general notion of participation, not simply for the case where
participation was set membership.}}

Russell's paradox drove set theory out of Frege's Paradise. As set theory was
reconstructed to escape the paradoxes, the set-universal for a property was
always non-self-participating. Thus set theory became not \textquotedblleft%
\textit{the} theory of universals\textquotedblright\ but the theory of
non-self-participating universals.

The reformulation of set theory cleared the ground for a separate theory of
always-self-participating universals. That idea was realized in category
theory $\cite{ell:cu}$ by the objects having universal mapping properties. The
self-participating universal for a property (if it exists) is the paradigmatic
or archetypical example of the property. All instances of the property are
determined to have the property by a morphism \textquotedblleft
participating\textquotedblright\ in that paradigmatic instance (where the
universal \textquotedblleft participates\textquotedblright\ in itself by the
identity morphism). The logic of category theory is the immanent logic of
determination through morphisms---particularly through universal morphisms
(see example below).

The notion of determination-by-morphisms plays the conceptually primitive role
in category theory analogous to the primitive notion of collection in set
theory \cite{fef:cf}. The set universal for a property collects together the
instances of the property but does not have the property itself. The
category-universal for a property has the property itself and determines the
instances of the property by morphisms.

The intuitive notion of a concrete universal{\footnote{I am ignoring some
rather woolly notions of \textquotedblleft concrete
universal\textquotedblright\ that exist in the philosophical literature other
than this notion of the concrete universal as the paradigm or exemplary
instance of a property. Also since the adjectives \textquotedblleft
concrete\textquotedblright\ and \textquotedblleft abstract\textquotedblright%
\ are already freighted with other meanings in category theory, I will refer
to the concrete universals in that context as self-participating universals.
Such a universal is \textquotedblleft concrete\textquotedblright\ only in the
sense that it is one among the instances of the property. Certainly all the
examples we meet in category theory are otherwise quite abstract mathematical
entities.} }occurs in ordinary language (any archetypical or paradigmatic
reference such as the \textquotedblleft all-American boy\textquotedblright\ or
the \textquotedblleft perfect\textquotedblright\ example of something), in the
arts and literature (the old idea that great art uses a concrete instance to
universally exemplify certain human conditions, e.g., Shakespeare's
\textit{Romeo and Juliet} exemplifies romantic tragedies), and in philosophy
(the pure example of F-ness with no imperfections, no junk, and no noise; only
those attributes necessary for F-ness). Some properties are even defined by
means of the concrete universal such as the property of being Lincolnesqe.
Abraham Lincoln is the concrete universal for the property and all other
persons with the property have it by virtue of resembling the concrete
universal. The vague intuitive notion of a concrete or self-participating
universal becomes quite precise in the universal mapping properties of CT.

Perhaps the breakthrough was MacLane's characterization of the direct product
$X\times Y$ of (say) two sets $X$ and $Y$ by a UMP. The property in question
is being \textquotedblleft a pair of maps, one to $X$ and one to $Y$, with a
common domain.\textquotedblright\ The self-participating universal is the
universal object $X\times Y$ and the pair of projections $p_{X}:X\times
Y\rightarrow X$ and $p_{Y}:X\times Y\rightarrow Y$. Given any other pair with
the property, $f:W\rightarrow X$ and $g:W\rightarrow Y$, there is a unique
factor map $h:W\rightarrow X\times Y$ such that $p_{X}h=f$ and $p_{Y}h=g$. The
pair $(f,g)$ are said to \textit{uniquely factor through} the projections
$(p_{X},p_{Y})$ by the factor map $h$. Thus a pair $(f,g)$ has the property if
and only if it participates in (uniquely factors through) the
self-participating universal $(p_{X},p_{Y})$. The UMP description of the
product characterizes it up to isomorphism but it does not show existence.

I will give a conceptual description of its construction that will be of use
later. Maps carry determination from one object to another. In this case, we
are considering determinations in both $X$ and $Y$ by some common domain set.
At the most `atomic' level, the determiner would be a one-point set which
would pick out `determinees', a point $x$ in $X$ and a point $y$ in $Y$. Thus
the most atomic determinees are the ordered pairs $(x,y)$ of elements from $X$
and $Y$. What would be the most universal determiner of $(x,y)$? Universality
is achieved by making the conceptual move of reconceptualizing the atoms
$(x,y)$ from being just a determinee to being its own determiner, i.e.,
self-determination. Thus the universal determiner of $X$ and $Y$ is the set of
all the ordered pairs $(x,y)$ from $X$ and $Y$, denoted $X\times Y $. It has
the UMP for the product and all other sets with that UMP would be isomorphic
to it.

\textquotedblleft Determination through self-determining
universals\textquotedblright\ sounds off-puttingly `philosophical' but it
nevertheless provides the abstract conceptual framework for the theory of
universals and adjoint functors, and thus for the view of category theory
presented here. For instance, a determinative relation has both a determiner
or sending end and a determinee or receiving end so there should be a pair of
universals, one at each end of the relationship, and that is precisely the
conceptual source of the pair of universals in an adjunction. But that is
getting ahead of the story.

Consider the dual case of the coproduct or disjoint union of sets. What is the
self-participating universal for two sets $X$ and $Y$ to determine a common
codomain set by a map from $X$ and a map from $Y$? At the most atomic level,
what is needed to determine a single point in the codomain set? A single point
$x$ in $X$ or a single point $y$ in $Y$ mapped to that single point would be
sufficient to determine it. Hence those single points from $X$ or $Y$ are the
atomic determiners or \textquotedblleft germs\textquotedblright\ to determine
a single point. What would be the most universal determinee of $x$ or of $y$?
Universality is achieved by making the conceptual move of reconceptualizing
$x$ and $y$ from being determiners to being their own determinees, i.e.,
self-determination. Thus the universal determinee of $X$ and $Y$ is the set of
all the $x$ from $X$ and the $y$ from $Y$, denoted $X+Y$. Since the two maps
from $X$ and $Y$ can be defined separately, any elements that might be common
to the two sets (i.e., in the intersection of the sets) are separate
determiners and thus would be separate determinees in $X+Y$. Thus it is not
the union but the disjoint union of $X$ and $Y$.

The universal instance of a pair of maps from $X$ and $Y$ to a common codomain
set is the pair of injections $i_{X}:X\rightarrow X+Y$ and $i_{Y}:Y\rightarrow
X+Y$ whereby each determiner $x$ in $X$ or $y$ in $Y$ determine themselves as
determinees in $X+Y$. For any other pair of maps $f:X\rightarrow Z$ and
$g:Y\rightarrow Z$ from $X$ and $Y$ to a common codomain set $Z$, there is a
unique factor map $h:X+Y\rightarrow Z$ such that $hi_{X}=f$ and $hi_{Y}=g$.
All determination of a common codomain set by $X$ and $Y$ flows or factors
through its own self-determination $i_{X}:X\rightarrow X+Y$ and $i_{Y}%
:Y\rightarrow X+Y$.

The product (or coproduct) universal could be used to illustrate the point
about the logic of category theory being the immanent logic of determination
through morphisms and particularly through universals. Consider an inference
in conventional logic: all roses are beautiful, $r$ is a rose, and therefore
$r$ is beautiful. This could be formulated in the \textit{language} of
category theory as follows: there is an inclusion of the set of roses $R$ in
the set of beautiful things $B$, i.e., $R\hookrightarrow B$, $r$ is an element
of $R$, i.e., $1\overset{r}{\rightarrow}R$, so by composition $r$ is an
element of $B$, i.e., $1\overset{r}{\rightarrow}R\hookrightarrow B=$
$1\overset{r}{\rightarrow}B$. This type of reformulation can be highly
generalized in category theory but my point is that this is only using the
powerful concepts of category theory to formulate and generalize conventional logic.

The immanent logic of category theory is different. Speaking in `philosophical
mode' for the moment, suppose there are two (self-participating or `concrete')
universals \textquotedblleft The Rose\textquotedblright\ $U_{R}$ and
\textquotedblleft The Beautiful" $U_{B}$ and that the The Rose participates in
The Beautiful by a morphism $U_{R}\rightarrow U_{B}$ which determines that The
Rose is beautiful. A particular entity $r$ participates in The Rose by a
morphism $r\rightarrow U_{R}$ which determines that $r$ is a rose. By
composition, the particular rose $r$ participates in The Beautiful by a
morphism $r\rightarrow U_{R}\rightarrow U_{B}=r\rightarrow U_{B}$ which
determines that the rose is beautiful. The logical inference is immanent in
the determinative morphisms.

For a categorical version, replace the property of being beautiful by the
property for which the product projections $(p_{X},p_{Y})$ were the universal,
i.e., the property of being \textquotedblleft a pair of maps $(f,g)$, one to
$X$ and one to $Y$, with a common domain, e.g., $f:W\rightarrow X$ and
$g:W\rightarrow Y$.\textquotedblright\ Suppose we have another map
$h:X\rightarrow Y$. \ Replace the property of being a rose with the property
of being \textquotedblleft a pair of maps $(f,g)$, one to $X$ and one to $Y$,
with a common domain such that $hf=g$, i.e., $W\overset{f}{\rightarrow
}X\overset{h}{\rightarrow}Y=W\overset{g}{\rightarrow}Y$. The universal for
that property (\textquotedblleft The Rose\textquotedblright) is given by a
pair of projections $(\pi_{X},\pi_{Y})$ with a common domain $Graph(h)$, the
\textit{graph} of $h$. There is a unique map $Graph(h)\overset{p}{\rightarrow
}X\times Y$ so that the graph\ participates in the product, i.e.,
$Graph(h)\overset{p}{\rightarrow}X\times Y\overset{p_{X}}{\rightarrow
}X=Graph(h)\overset{\pi_{X}}{\rightarrow}X$ and similarly for the other
projection. Now consider a particular \textquotedblleft rose\textquotedblright%
, namely a pair of maps $r_{X}:W\rightarrow X$ and $r_{Y}:W\rightarrow Y$ that
commute with $h$ ($hr_{X}=r_{Y}$) and thus participate in the graph, i.e.,
there is a unique map $r:W\rightarrow Graph(h)$ such that $W\overset
{r}{\rightarrow}Graph(h)\overset{\pi_{X}}{\longrightarrow}X=W\overset{r_{X}%
}{\longrightarrow}X$ and $W\overset{r}{\rightarrow}Graph(h)\overset{\pi_{Y}%
}{\longrightarrow}Y=W\overset{r_{Y}}{\longrightarrow}Y$. Then by composition,
the particular pair of maps $(r_{X},r_{Y})$ also participates in the product,
i.e., $pr$ is the unique map such that $p_{X}pr=r_{X}$ and $p_{Y}pr=r_{Y}$.
Not only does having the property imply participation in the universal, any
entity that participates in the universal is thereby forced to have the
property. For instance, because $(r_{X},r_{Y})$ participates in the product,
it has to be a pair of maps to $X$ and to $Y$ with a common domain (the
property represented by the product).

This example of the immanent logic of category theory, where maps play a
determinative role and the determination is through universals, could be
contrasted to formulating and generalizing the set treatment of the inference
in the language of categories. Let $P$ be the set of pairs of maps to $X$ and
$Y$ with a common domain, let $G$ be the set of pairs of such maps to $X$ and
$Y$ that also commute with $h$, and let $(r_{X},r_{Y})$ be a particular pair
of such maps. Then $1\overset{(r_{X},r_{Y})}{\longrightarrow}G\hookrightarrow
P=1\rightarrow P$. That uses the language of category theory to formulate the
set treatment with the set-universals for those properties and it ignores the
category-universals for those properties.

Our purpose here is not to further develop the (immanent) logic of category
theory but to show how adjoint functors fit within that general framework of
determination through universals.

\section{Definition and Directionality of Adjoints}

There are many equivalent definitions of adjoint functors (see MacLane
$\cite{mac:cwm}$), but the most `official' one seems to be the one using a
natural isomorphism of hom-sets. Let $\mathbf{X}$ and $\mathbf{A}$ be
categories and $F:\mathbf{X\rightarrow A}$ and $G:\mathbf{A\rightarrow X}$
functors between them. Then $F$ and $G$ are said to be a pair of
\textit{adjoint functors} or an \textit{adjunction}, written $F\dashv G$, if
for any $x$ in $\mathbf{X}$ and $a$ in $\mathbf{A}$, there is an isomorphism
$\phi$ natural in $x$ and in $a$:%

\[
{\phi}_{x,a}:Hom_{\mathbf{A}}(Fx,a)\;{\cong\;}Hom_{\mathbf{X}}(x,Ga).
\]

\noindent With this standard way of writing the isomorphism of hom-sets, the
functor $F$ on the left is called the \textit{left adjoint} and the functor
$G$ on the right is the \textit{right adjoint}. Maps associated with each
other by the adjunction isomorphism ("adjoint transposes" of one another) are
indicated by an asterisk so if $g:Fx\rightarrow a$ then $g^{\ast}:x\rightarrow
Ga$ is the associated map ${\phi}_{x,a}(g)=g^{\ast}$ and similarly if
$f:x\rightarrow Ga$ then $\phi_{x,a}^{-1}(f)=f^{\ast}:Fx\rightarrow a$ is the
associated map.

In much of the literature, adjoints are presented in a seemingly symmetrical
fashion so that there appears to be no directionality of the adjoints between
the categories $\mathbf{X}$ and $\mathbf{A}$. But there is a directionality
and it is important in understanding adjoints. Both the maps that appear in
the adjunction isomorphism, $Fx\rightarrow a$ and $x\rightarrow Ga$, go from
the \textquotedblleft$x$-thing\textquotedblright\ (i.e., either $x$ or the
image $Fx$) to the \textquotedblleft$a$-thing\textquotedblright\ (either the
image $Ga$ or $a$ itself), so we see a direction emerging from $\mathbf{X}$ to
$\mathbf{A}$. That direction of an adjunction is the direction of the left
adjoint (which goes from $\mathbf{X}$ to $\mathbf{A}$). Then $\mathbf{X}$
might called the \textit{sending} category and $\mathbf{A}$ the
\textit{receiving} category.\footnote{Sometimes adjunctions are written with
this direction as in the notation $\langle F,G,{\phi\rangle}:\mathbf{X}%
{\rightharpoonup}$ $\mathbf{A}$ (MacLane $\cite[p. 78]{mac:cwm}$). This also
allows the composition of adjoints to be defined in a straightforward manner
(MacLane $\cite[p. 101]{mac:cwm}$).}

Bidirectionality of determination through adjoints occurs when a functor has
both a left and right adjoint. \ We will later see such an example of
bidirectionality of determination between sets and diagram functors where the
limit and colimit functors are respectively right and left adjoints to the
same constant functor.\footnote{Also any adjunction can be restricted to
subcategories where the unit and counit are isomorphisms so that each adjoint
is both a left and right adjoint of the other---and thus determination is
bidirectional---on those subcategories (see \cite{lam:her} or \cite{bel:top}%
).}

In the theory of adjoints presented here, the directionality of adjoints
results from being representations of heteromorphisms\ which have that
directionality. \ Such morphisms can exhibited in concrete examples of
adjoints (see the later examples). To abstractly define chimera morphisms\ or
heteromorphisms that work for all adjunctions, we turn to the presentation of
adjoints using adjunctive squares.

\section{Adjunctive Squares}

\subsection{Embedding Adjunctions in a Product Category}

Our approach to a theory of adjoints uses a certain \textquotedblleft
adjunctive square\textquotedblright\ diagram that is in the product category
$\mathbf{X}\times\mathbf{A}$ associated with an adjunction $F:\mathbf{X}%
\rightleftarrows\mathbf{A}:G$. With each object $x$ in the category
$\mathbf{X}$, we associate the element $\widehat{x}=(x,Fx)$ in the product
category $\mathbf{X}\times\mathbf{A}$ so that $Ga$ would have associated with
it $\widehat{Ga}=(Ga,FGa)$.\ With each morphism in $\mathbf{X}$ with the form
$h:x^{\prime}\rightarrow x$, we associate the morphism $\widehat
{h}=(h,Fh):\widehat{x^{\prime}}=(x^{\prime},Fx^{\prime})\rightarrow\widehat
{x}=(x,Fx)$ in the product category $\mathbf{X}\times\mathbf{A}$ (maps compose
and diagrams commute component-wise). Thus the mapping of $x$ to $(x,Fx)$
extends to an embedding $(1_{\mathbf{X}},F):\mathbf{X}\rightarrow
\mathbf{X}\times\mathbf{A}$ whose image $\widehat{\mathbf{X}}$ is isomorphic
with $\mathbf{X}$.

With each object $a$ in the category $\mathbf{A}$, we associate the element
$\widehat{a}=(Ga,a)$ in the product category $\mathbf{X}\times\mathbf{A}$ so
that $Fx$ would have associated with it $\widehat{Fx}=(GFx,Fx)$. With each
morphism in $\mathbf{A}$ with the form $k:a\rightarrow a^{\prime},$we
associate the morphism $\widehat{k}=(Gk,k):(Ga,a)\rightarrow(Ga^{\prime
},a^{\prime})$ in the product category $\mathbf{X}\times\mathbf{A}$. The
mapping of $a$ to $(Ga,a)$ extends to an embedding $(G,1_{\mathbf{A}%
}):\mathbf{A}\rightarrow\mathbf{X}\times\mathbf{A}$ whose image $\widehat
{\mathbf{A}}$ is isomorphic to $\mathbf{A}$.

The mapping of $x$ to $(GFx,Fx)$ in $\widehat{\mathbf{A}}$ extends to the
functor $(GF,F):\mathbf{X}\rightarrow\mathbf{X}\times\mathbf{A}$. Given that
$F$ and $G$ are adjoints, the unit natural transformation $\eta:1_{\mathbf{X}%
}\rightarrow GF$ of the adjunction and the identity natural transformation
$1_{F}$ give a natural transformation $(\eta,1_{F}):(1_{\mathbf{X}%
},F)\rightarrow(GF,F)$ with the component at $x$ being $(\eta_{x}%
,1_{Fx}):(x,Fx)\rightarrow(GFx,Fx)$.

The mapping of $a$ to $(Ga,FGa)$ in $\widehat{\mathbf{X}}$ extends to a
functor $(G,FG):\mathbf{A}\rightarrow\mathbf{X}\times\mathbf{A}$. The counit
natural transformation $\varepsilon$ of the adjunction and the identity
natural transformation $1_{G}$ give a natural transformation $(1_{G}%
,\varepsilon):(G,FG)\rightarrow(G,1_{\mathbf{A}})$ with the component at $a$
being $(1_{G},\varepsilon_{a}):(Ga,FGa)\rightarrow(Ga,a)$.

The $(F,G)$ twist functor, which carries $(x,a)$ to $(Ga,Fx)$, is an
endo-functor on $\mathbf{X}\times\mathbf{A}$ which carries $\widehat
{\mathbf{X}}$ to $\widehat{\mathbf{A}}$ and $\widehat{\mathbf{A}}$ to
$\widehat{\mathbf{X}}$ to reproduce the adjunction between those two subcategories.

These various parts can then be collected together in the adjunctive square diagram.

\begin{center}
$%
\begin{array}
[c]{rcccl}
& (x,Fx) & \overset{(f,Ff)}{\rightarrow} & (Ga,FGa) & \\
(\eta_{x},1_{Fx}) & \downarrow &  & \downarrow & (1_{Ga},\varepsilon_{a})\\
& (GFx,Fx) & \overset{(Gg,g)}{\longrightarrow} & (Ga,a) &
\end{array}
\medskip$

Adjunctive Square Diagram
\end{center}

Any diagrams of this form for some given $f:x\rightarrow Ga$ or for some given
$g:Fx\rightarrow a$ will be called \textit{adjunctive squares}. The purpose of
the adjunctive square diagram is to conveniently represent the properties of
an adjunction in the format of commutative squares. \ The map on the top is in
$\widehat{\mathbf{X}}$ (the \textquotedblleft top category\textquotedblright)
and the map on the bottom is in $\widehat{\mathbf{A}}$ (the \textquotedblleft
bottom category\textquotedblright) and the vertical maps as well as the main
diagonal $(f,g)$ in a commutative adjunctive square are morphisms from
$\widehat{\mathbf{X}}$-objects to $\widehat{\mathbf{A}}$-objects. The $(F,G)$
twist functor carries the main diagonal map $(f,g)$ to the (southwest to
northeast) anti-diagonal map $(Gg,Ff)$.

Given $f:x\rightarrow Ga$, the rest of the diagram is determined by the
requirement that the square commutes. Commutativity in the second component
uniquely determines that $g=g1_{Fx}=\varepsilon_{a}Ff$ so $g=f^{\ast
}=\varepsilon_{a}Ff$ \ is the map associated with $f$ in the adjunction
isomorphism. \ Commutativity in the first component is the universal mapping
property factorization of the given $f:x\rightarrow Ga$ through the unit
$x\overset{\eta_{x}}{\longrightarrow}GFx\overset{Gf^{\ast}}{\longrightarrow
}Ga=x\overset{f}{\longrightarrow}Ga$. Similarly, if we were given
$g:Fx\rightarrow a,$ then commutativity in the first component implies that
$f=1_{Ga}f=Gg\eta_{x}=g^{\ast}$. And commutativity in the second component is
the UMP factorization of $g:Fx\rightarrow a$ through the counit $Fx\overset
{Fg^{\ast}}{\longrightarrow}FGa\overset{\varepsilon_{a}}{\longrightarrow
}a=Fx\overset{g}{\longrightarrow}a$.

The adjunctive square diagram also brings out the directionality of the
adjunction. In a commutative adjunctive square, the main diagonal map, which
goes from $\widehat{x}=(x,Fx)$ in $\widehat{\mathbf{X}}$ to $\widehat
{a}=(Ga,a)$ in $\widehat{\mathbf{A}}$, is $(f,g)$ where $g=f^{\ast}$ and
$f=g^{\ast}$. Each $(x,Fx)$-to-$(Ga,a)$ determination crosses one of the
\textquotedblleft bridges\textquotedblright\ represented by the natural
transformations $(\eta,1_{F})$ or $(1_{G},\varepsilon)$. Suppose that a
determination went from $(x,Fx)$ to $(Ga,a)$ as follows: $(x,Fx)\rightarrow
(x^{\prime},Fx^{\prime})\rightarrow(GFx^{\prime},Fx^{\prime})\rightarrow
(Ga,a)$.

\begin{center}
$%
\begin{array}
[c]{ccccc}%
(x,Fx) & \longrightarrow & (x^{\prime},Fx^{\prime}) &  & \\
\downarrow &  & \downarrow &  & \\
(GFx,Fx) & \longrightarrow & (GFx^{\prime},Fx^{\prime}) & \longrightarrow &
(Ga,a)
\end{array}
$
\end{center}

\noindent By the commutativity of the square from the naturality of
$(\eta,1_{F})$, the determination $(x,Fx)$ to $(Ga,a)$ crossing the bridge at
$(x,Fx)\rightarrow(GFx,Fx)$ is the same. Hence any determination from $(x,Fx)
$ and crossing a $(\eta,1_{F})$-bridge can be taken as crossing the bridge at
$(x,Fx)$.

Suppose a determination from $(x,Fx)$ to $(Ga,a)$ crossed a $(1_{G}%
,\varepsilon)$-bridge: $(x,Fx)\rightarrow(Ga^{\prime},FGa^{\prime}%
)\rightarrow(Ga^{\prime},a^{\prime})\rightarrow(Ga,a)$.

\begin{center}
$%
\begin{array}
[c]{ccccc}%
(x,Fx) & \longrightarrow & (Ga^{\prime},FGa^{\prime}) & \longrightarrow &
(Ga,FGa)\\
&  & \downarrow &  & \downarrow\\
&  & (Ga^{\prime},a^{\prime}) & \longrightarrow & (Ga,a)
\end{array}
$
\end{center}

\noindent By the commutativity of the square from the naturality of
$(1_{G},\varepsilon)$, the determination crossing the bridge at
$(Ga,FGa)\rightarrow(Ga,a)$ is the same. Hence any determination to $(Ga,a)$
and crossing a $(1_{G},\varepsilon)$-bridge can be taken as crossing the
bridge at $(Ga,a)$. Thus the adjunctive square represents the general case of
possible $(x,Fx)$-to-$(Ga,a)$ determinations using functors $F$ and $G$.

\subsection{Factorization Systems of Maps}

Consider the generic form for an adjunctive square diagram where we are
assuming that $F$ and $G$ are adjoint.

\begin{center}
$%
\begin{array}
[c]{rcccl}
& (x,Fx) & \overset{(f,Ff)}{\rightarrow} & (Ga,FGa) & \\
(\eta_{x},1_{Fx}) & \downarrow & \nearrow & \downarrow & (1_{Ga}%
,\varepsilon_{a})\\
& (GFx,Fx) & \overset{(Gg,g)}{\longrightarrow} & (Ga,a) &
\end{array}
$
\end{center}

Given any $f:x\rightarrow Ga$ in the first component on top, there is a unique
$f^{\ast}:Fx\rightarrow a$ such that the anti-diagonal map $(Gf^{\ast},Ff)$
factors $(f,Ff)$ through the universal $(\eta_{x},1_{Fx})$, i.e., such that
the upper triangle commutes. And that anti-diagonal map composes with
$(1_{Ga},\varepsilon_{a})$ to determine the bottom map

\begin{center}
$(GFx,Fx)\overset{(Gf^{\ast},Ff)}{\longrightarrow}(Ga,FGa)\overset
{(1_{Ga},\varepsilon_{a})}{\longrightarrow}(Ga,a)=(GFx,Fx)\overset{(Gf^{\ast
},f^{\ast})}{\longrightarrow}(Ga,a)$
\end{center}

\noindent that makes the lower triangle and thus the square commutes. Or
starting with $g=f^{\ast}:Fx\rightarrow a$ in the second component on the
bottom, there is a unique $g^{\ast}=f:x\rightarrow Ga$ such that the
anti-diagonal map $(Gf^{\ast},Ff)$ factors $(Gg,g)=(Gf^{\ast},f^{\ast})$
through the universal $(1_{Ga},\varepsilon_{a})$, i.e., such that the lower
triangle commutes. \ And that anti-diagonal map composes with $(\eta
_{x},1_{Fx})$ to determine the upper map

\begin{center}
$(x,Fx)\overset{(\eta_{x},1_{Fx})}{\longrightarrow}(GFx,Fx)\overset{(Gf^{\ast
},Ff)}{\longrightarrow}(Ga,FGa)=(x,Fx)\overset{(f,Ff)}{\longrightarrow
}(Ga,FGa)$
\end{center}

\noindent that makes the upper triangle and thus the square commutes.

Adjointness is closely related to the notion of a factorization system for
orthogonal sets of maps such as epis and monos (see $\cite{Freyd:ccfI}$ or
$\cite{tay:pfm}$). The motivating example for factorization systems of maps
was the example of epimorphisms and monomorphisms. The unit and counit in any
adjunction have a closely related property. For instance, suppose there are
morphisms $g,\ g^{\prime}:Fx\rightarrow a$ in the bottom category $\mathbf{A}$
such that $x\overset{\eta_{x}}{\longrightarrow}GFx\overset{Gg}{\longrightarrow
}Ga=x\overset{\eta_{x}}{\longrightarrow}GFx\overset{Gg^{\prime}}%
{\longrightarrow}Ga$ holds in the top category $\mathbf{X}$. By the unique
factorization of any morphism $f:x\rightarrow Ga$ through the unit $\eta_{x}$
it follows that $g=g^{\prime}$. Since the uniqueness of the morphisms is in
$\mathbf{A}$ while the $G$-image maps $Gg$ and $Gg^{\prime}$ post-composed
with $\eta_{x}$ are in $\mathbf{X}$, this is not the same as $\eta_{x}$ being
epi. But it is a closely related property of being \textquotedblleft epi with
respect to $G$-images of $Fx\rightarrow a$ morphisms.\textquotedblright%
\ Similarly, the counit $\varepsilon_{a}$ is \textquotedblleft mono with
respect to $F$-images of $x\rightarrow Ga$ morphisms.\textquotedblright\ Thus
it is not surprising that the maps $(\eta_{x},1_{Fx})$ and $(1_{Ga}%
,\varepsilon_{a})$ have an anti-diagonal factorization in the adjunctive
square diagrams in a manner analogous to epis and monos in general commutative squares.

\section{Adjoints = Birepresentations of Het-bifunctors}

\subsection{Chimera Morphisms and Het-bifunctors}

We argue that an adjunction has to do with morphisms between objects that are
in general in different categories, e.g., from an object $x$ in $\mathbf{X}$
to an object $a$ in $\mathbf{A}$. These \textquotedblleft
heteromorphisms\textquotedblright\ (in contrast to homomorphisms) are like
mongrels or chimeras that do not fit into either of the two categories. For
instance, in the context of the free-group-underlying-set adjunction, we might
consider any mapping $x\overset{c}{\Rightarrow}a$ whose tail is a set and head
is a group. Some adjunctions, such as the free-underlying adjunctions, have
this sort of concrete realization of the cross-category determination which is
then represented by the two adjoint transposes $x\overset{f(c)}{\rightarrow
}Ga$ and $Fx\overset{g(c)}{\rightarrow}a$. Since the cross-category
heteromorphisms are not morphisms in either of the categories, what can we say
about them?

The one thing we can reasonably say is that chimera morphisms can be
precomposed or postcomposed with morphisms within the categories (i.e.,
intra-category morphisms) to obtain other chimera morphisms.\footnote{The
chimera genes are dominant in these mongrel matings. While mules cannot mate
with mules, it is `as if' mules could mate with either horses or donkeys to
produce other mules.} This is easily formalized. Suppose we have
heteromorphisms $x\overset{c}{\Rightarrow}a$ as in the case of sets and groups
from objects in $\mathbf{X}$ to objects in $\mathbf{A}$ (another example
analyzed below would be \textquotedblleft cones\textquotedblright\ which can
be seen as chimera morphisms from sets to diagram functors). Let
$Het(x,a)=\{x\overset{c}{\Rightarrow}a\}$ be the set of heteromorphisms from
$x$ to $a$. For any $\mathbf{A}$-morphism $k:a\rightarrow a^{\prime}$ and any
chimera morphism $x\overset{c}{\Rightarrow}a$, intuitively there is a
composite chimera morphism $x\overset{c}{\Rightarrow}a\overset{k}{\rightarrow
}a^{\prime}=x\overset{kc}{\Rightarrow}a^{\prime}$, i.e., $k$ induces a map
$Het(x,k):Het(x,a)\rightarrow Het(x,a^{\prime})$. For any $\mathbf{X}%
$-morphism $h:x^{\prime}\rightarrow x$ and chimera morphism $x\overset
{c}{\Rightarrow}a$, intuitively there is the composite chimera morphism
$x^{\prime}\overset{h}{\rightarrow}x\overset{c}{\Rightarrow}a=x^{\prime
}\overset{ch}{\Rightarrow}a$, i.e., $h$ induces a map
$Het(h,a):Het(x,a)\rightarrow Het(x^{\prime},a)$ (note the reversal of
direction). Taking the sets-to-groups example to guide intuitions, the induced
maps would respect identity and composite morphisms in each category.
Moreover, composition is associative in the sense that $(kc)h=k(ch) $. This
means that the assignments of sets of chimera morphisms $Het(x,a)=\{x\overset
{c}{\Rightarrow}a\}$ and the induced maps between them constitute a bifunctor
$Het:\mathbf{X}^{op}\times\mathbf{A}\rightarrow\mathbf{Set}$ (contravariant in
the first variable and covariant in the second).\footnote{One conventional
treatment of chimera morphisms such as cones `misses' the bifunctor
formulation by treating the cones as objects in a category (rather than as
morphisms between categories). The compositions are then viewed as defining
morphisms between these objects. For instance, $h:x^{\prime}\rightarrow x$
would define a morphism from a cone $ch:x^{\prime}\Rightarrow a$ to the cone
$c:x\Rightarrow a$. The terminal object in this category would be the limit
cone---and dually for cocones. These are slices of the chimera comma category
defined in the next section.} The composition properties we would axiomatize
for \textquotedblleft chimera morphisms\textquotedblright\ or
\textquotedblleft heteromorphisms\textquotedblright\ from $\mathbf{X}$-objects
to $\mathbf{A}$-objects are precisely those of the elements in the values of
such a bifunctor $Het:\mathbf{X}^{op}\times\mathbf{A}\rightarrow\mathbf{Set}$.

An adjunction is not simply about heteromorphisms from $\mathbf{X}$ to
$\mathbf{A}$; it is about such determinations through universals. In other
words, an adjunction arises from a het-bifunctor $Het:\mathbf{X}^{op}%
\times\mathbf{A}\rightarrow\mathbf{Set}$ that is \textquotedblleft
birepresentable\textquotedblright\ in the sense of being representable on both
the left and right.

Given any bifunctor $Het:\mathbf{X}^{op}\times\mathbf{A}\rightarrow
\mathbf{Set}$, it is \textit{representable on the left} if for each
$\mathbf{X}$-object $x$, there is an $\mathbf{A}$-object $Fx$ that represents
the functor $Het(x,-)$, i.e., there is an isomorphism $\psi_{x,a}%
:Hom_{\mathbf{A}}(Fx,a)\cong Het(x,a)$ natural in $a$. For each $x$, let
$h_{x}$ be the image of the identity on $Fx$, i.e., $\psi_{x,Fx}(1_{Fx}%
)=h_{x}\in Het(x,Fx) $. We first show that $h_{x}$ is a universal element for
the functor $Het(x,-)$ and then use that to complete the construction of $F$
as a functor. For any $c\in Het(x,a)$, let $g(c)=\psi_{x,a}^{-1}%
(c):Fx\rightarrow a$. Then naturality in $a$ means that the following diagram commutes.

\begin{center}
$%
\begin{array}
[c]{ccccc}
& Hom_{\mathbf{A}}(Fx,Fx) & \cong & Het(x,Fx) & \\
_{Hom(Fx,g(c))} & \downarrow &  & \downarrow & _{Het(x,g(c))}\\
& Hom_{\mathbf{A}}(Fx,a) & \cong & Het(x,a) &
\end{array}
$
\end{center}

\noindent Chasing $1_{Fx}$ around the diagram yields that $c=Het(x,g(c))(h_{x}%
)$ which can be written as $c=g(c)h_{x}$. Since the horizontal maps are
isomorphisms, $g(c)$ is the unique map $g:Fx\rightarrow a $ such that
$c=gh_{x}$. Then $(Fx,h_{x})$ is a \textit{universal element} (in MacLane's
sense \cite[p. 57]{mac:cwm}) for the functor $Het(x,-)$ or equivalently
$1\overset{h_{x}}{\longrightarrow}Het(x,Fx)$ is a \textit{universal arrow}
\cite[p. 58]{mac:cwm} from $1$ (the one point set) to $Het(x,-)$. Then for any
$\mathbf{X}$-morphism $j:x\rightarrow x^{\prime}$, $Fj:Fx\rightarrow
Fx^{\prime}$ is the unique $\mathbf{A}$-morphism such that $Het(x,Fj)$ fills
in the right vertical arrow in the following diagram.

\begin{center}
$%
\begin{array}
[c]{rcccc}
& 1 & \overset{h_{x}}{\longrightarrow} & Het(x,Fx) & \\
_{h_{x^{\prime}}} & \downarrow &  & \downarrow & _{Het(x,Fj)}\\
& Het(x^{\prime},Fx^{\prime}) & \overset{Het(j,Fx^{\prime})}{\longrightarrow}
& Het(x,Fx^{\prime}) &
\end{array}
$
\end{center}

\noindent It is easily checked that such a definition of $Fj:Fx\rightarrow
Fx^{\prime}$ preserves identities and composition using the functoriality of
$Het(x,-)$ so we have a functor $F:X\rightarrow A$. It is a further standard
result that the isomorphism is also natural in $x$ (e.g., \cite[p.
81]{mac:cwm} or the "parameter theorem" \cite[p. 525]{macbirk:alg}).

Given a bifunctor $Het:\mathbf{X}^{op}\times\mathbf{A}\rightarrow\mathbf{Set}%
$, it is \textit{representable on the right} if for each $\mathbf{A}$-object
$a$, there is an $\mathbf{X}$-object $Ga$ that represents the functor
$Het(-,a)$, i.e., there is an isomorphism $\varphi_{x,a}:Het(x,a)\cong
Hom_{X}(x,Ga)$ natural in $x$. For each $a$, let $e_{a}$ be the inverse image
of the identity on $Ga$, i.e., $\varphi_{Ga,a}^{-1}(1_{Ga})=e_{a}\in
Het(Ga,a)$. For any $c\in Het(x,a)$, let $f(c)=\varphi_{x,a}(c):x\rightarrow
Ga$. Then naturality in $x$ means that the following diagram commutes.

\begin{center}
$%
\begin{array}
[c]{ccccc}
& Het(Ga,a) & \cong & Hom_{\mathbf{X}}(Ga,Ga) & \\
_{Het(f(c),a)} & \downarrow &  & \downarrow & _{Hom(f(c),Ga)}\\
& Het(x,a) & \cong & Hom_{\mathbf{X}}(x,Ga) &
\end{array}
$
\end{center}

\noindent Chasing $1_{Ga}$ around the diagram yields that $c=Het(f(c),a)(e_{a}%
)=e_{a}f(c)$ so $(Ga,e_{a})$ is a universal element for the functor $Het(-,a)$
and that $1\overset{e_{a}}{\longrightarrow}Het(Ga,a)$ is a universal arrow
from $1$ to $Het(-,a)$. Then for any $\mathbf{A}$-morphism $k:a^{\prime
}\rightarrow a$, $Gk:Ga^{\prime}\rightarrow Ga$ is the unique $\mathbf{X}%
$-morphism such that $Het(Gk,a)$ fills in the right vertical arrow in the
following diagram.

\begin{center}
$%
\begin{array}
[c]{ccccc}
& 1 & \overset{e_{a}}{\longrightarrow} & Het(Ga,a) & \\
_{e_{a^{\prime}}} & \downarrow &  & \downarrow & _{Het(Gk,a)}\\
& Het(Ga^{\prime},a^{\prime}) & \overset{Het(Ga^{\prime},k)}{\longrightarrow}
& Het(Ga^{\prime},a) &
\end{array}
$
\end{center}

\noindent In a similar manner, it is easily checked that the functoriality of
$G$ follows from the functoriality of $Het(-,a)$. Thus we have a functor
$G:\mathbf{A}\rightarrow\mathbf{X}$ such that $Ga$ represents the functor
$Het(-,a)$, i.e., there is a natural isomorphism $\varphi_{x,a}:Het(x,a)\cong
Hom_{\mathbf{X}}(x,Ga)$ natural in $x$. And in a similar manner, it can be
shown that the isomorphism is natural in both variables. Thus given a
bifunctor $Het:\mathbf{X}^{op}\times\mathbf{A}\rightarrow\mathbf{Set}$
representable on both sides, we have an adjunction isomorphism:%
\[
Hom_{\mathbf{A}}(Fx,a)\cong Het(x,a)\cong Hom_{\mathbf{X}}(x,Ga).
\]

\noindent The morphisms mapped to one another by these isomorphisms are
\textit{adjoint correlates} of one another. When the two representations are
thus combined, the universal element $h_{x}\in Het(x,Fx)$ induced by $1_{Fx}$
is the \textit{chimera unit} and is the adjoint correlate of the ordinary unit
$\eta_{x}:x\rightarrow GFx$. The universal element $e_{a}\in Het(Ga,a)$
induced by $1_{Ga}$ is the \textit{chimera counit} and is the adjoint
correlate of the ordinary counit $\varepsilon_{a}:FGa\rightarrow a$. The two
factorizations $g(c)h_{x}=c=e_{a}f(c)$ combine to give what we will later call
the \textquotedblleft chimera adjunctive square\textquotedblright\ with $c$ as
the main diagonal.

\subsection{Comma Categories and Bifunctors}

The above treatment uses the `official' hom-set definition of an adjunction.
\ It might be useful to briefly restate the ideas using William Lawvere's
comma category definition of an adjunction in his famous 1963 thesis
\cite{law:thesis}. Given three categories $\mathbf{A}$, $\mathbf{B}$, and
$\mathbf{C}$ and two functors%

\[
A:\mathbf{A}\rightarrow\mathbf{C}\leftarrow\mathbf{B}:B
\]

\noindent the \textit{comma category} $(A,B)$ has as objects the morphism in
$\mathbf{C}$ from a value $A(a)$ of the functor $A:\mathbf{A}\rightarrow
\mathbf{C}$ to a value $B(b)$ of the functor $B:\mathbf{B}\rightarrow
\mathbf{C}$. A morphism from an object $A(a)\rightarrow B(b)$ to an object
$A(a^{\prime})\rightarrow B(b^{\prime})$ is an $\mathbf{A}$-morphism
$k:a\rightarrow a^{\prime}$ and a $\mathbf{B}$-morphism $h:b\rightarrow
b^{\prime}$ such that the following diagram commutes:

\begin{center}
$%
\begin{array}
[c]{ccccc}
& A(a) & \overset{A(k)}{\longrightarrow} & A(a^{\prime}) & \\
& \downarrow &  & \downarrow & \\
& B(b) & \overset{B(h)}{\longrightarrow} & B(b^{\prime}). &
\end{array}
$
\end{center}

\noindent There are two projection functors $\pi_{0}:(A,B)\rightarrow
\mathbf{A}$ which takes the object $A(a)\rightarrow B(b)$ to the object $a$ in
$\mathbf{A}$ and a morphism $(k,h)$ to $k$, and $\pi_{1}:(A,B)\rightarrow
\mathbf{B}$ with the analogous definition.

Then with the adjunctive setup $F:\mathbf{X}\rightleftarrows\mathbf{A}:G$,
there are two associated comma categories $(F,1_{\mathbf{A}})$ defined by
$F:\mathbf{X}\rightarrow\mathbf{A}\leftarrow\mathbf{A}:1_{\mathbf{A}}$, and
$(1_{\mathbf{X}},G)$ defined by $1_{\mathbf{X}}:\mathbf{X}\rightarrow
\mathbf{X}\leftarrow\mathbf{A}:G$. Then Lawvere's definition of an adjunction
(e.g., \cite[p. 84]{mac:cwm}, \cite[p. 92]{mclarty:ecet}, or \cite[p.
389]{tay:pfm}) is that $F$ is the left adjoint to $G$ if there is an
isomorphism $(F,1_{\mathbf{A}})\cong(1_{\mathbf{X}},G)$ over the projections
to the product category $\mathbf{X}\times\mathbf{A}$, i.e., such that the
following diagram commutes.

\begin{center}
$%
\begin{array}
[c]{ccccc}
& (F,1_{\mathbf{A}}) & \cong & (1_{\mathbf{X}},G) & \\
(\pi_{0},\pi_{1}) & \downarrow &  & \downarrow & (\pi_{0},\pi_{1})\\
& \mathbf{X}\times\mathbf{A} & = & \mathbf{X}\times\mathbf{A} &
\end{array}
$
\end{center}

To connect the comma category approach to the bifunctor treatment, we might
note that the data $A:\mathbf{A}\rightarrow\mathbf{C}\leftarrow\mathbf{B}:B $
used to define the comma category also defines a bifunctor $Hom_{\mathbf{C}%
}(A(-),B(-)):\mathbf{A}^{op}\times\mathbf{B}\rightarrow\mathbf{Set}$. But we
might start with any bifunctor $\mathcal{D}:\mathbf{X}^{op}\times
\mathbf{A}\rightarrow\mathbf{Set}$ and then mimic the definition of a comma
category to arrive at a category that we will denote $\mathcal{D}%
(1_{\mathbf{X}},1_{\mathbf{A}})$. The objects of $\mathcal{D}(1_{\mathbf{X}%
},1_{\mathbf{A}})$ are elements of any value $\mathcal{D}(x,a)$ of the
bifunctor which we will denote $x\overset{c}{\Rightarrow}a$. A morphism from
$x\overset{c}{\Rightarrow}a$ to $x^{\prime}\overset{c^{\prime}}{\Rightarrow
}a^{\prime} $ is given by an $\mathbf{X}$-morphism $j:x\rightarrow x^{\prime}$
and an $\mathbf{A}$-morphism $k:a\rightarrow a^{\prime}$ such that the
following diagram commutes:

\begin{center}
$%
\begin{array}
[c]{ccccc}
& x & \overset{j}{\longrightarrow} & x^{\prime} & \\
c & \Downarrow &  & \Downarrow & c^{\prime}\\
& a & \overset{k}{\longrightarrow} & a^{\prime} &
\end{array}
$
\end{center}

\noindent which means in terms of the bifunctor that the following diagram commutes:

\begin{center}
$%
\begin{array}
[c]{rcccc}
& 1 & \overset{c^{\prime}}{\longrightarrow} & \mathcal{D}(x^{\prime}%
,a^{\prime}) & \\
_{c} & \downarrow &  & \downarrow & _{\mathcal{D}(j,a^{\prime})}\\
& \mathcal{D}(x,a) & \overset{\mathcal{D}(x,k)}{\longrightarrow} &
\mathcal{D}(x,a^{\prime}). &
\end{array}
$
\end{center}

\noindent The projection functors to $\mathbf{X}$ and to $\mathbf{A}$ are
defined in the obvious manner.

There is a somewhat pedantic question of whether or not $\mathcal{D}%
(1_{\mathbf{X}},1_{\mathbf{A}})$ should be called a \textquotedblleft comma
category.\textquotedblright\ Technically it is not since its objects are
elements of an arbitrary bifunctor $\mathcal{D}:\mathbf{X}^{op}\times
\mathbf{A}\rightarrow\mathbf{Set}$, not morphisms within any category
$\mathbf{C}$. This point is not entirely pedantic for our purposes since when
the elements of $\mathcal{D}(x,a)$ are taken to be concrete chimera morphisms,
then they indeed are not morphisms within a category but are
(object-to-object) heteromorphisms between categories. For instance, there is
no comma category $(1_{\mathbf{X}},1_{\mathbf{A}})$ with identity functors on
both sides unless $\mathbf{X}$ and $\mathbf{A}$ are the same category.
However, $\mathcal{D}(1_{\mathbf{X}},1_{\mathbf{A}})$ will often be isomorphic
to a comma category (e.g., when $\mathcal{D}$ is representable on one side or
the other) and Lawvere only considered that comma categories were defined
\textquotedblleft up to isomorphism.\textquotedblright\ Moreover, I think that
ordinary usage by category theorists would call $\mathcal{D}(1_{\mathbf{X}%
},1_{\mathbf{A}})$ a \textquotedblleft comma category\textquotedblright\ even
though its objects are not necessarily morphisms within any category
$\mathbf{C}$, and I will follow that general usage. Thus $\mathcal{D}%
(1_{\mathbf{X}},1_{\mathbf{A}})$ is the \textit{comma category of the
bifunctor }$\mathcal{D}$. This broader notion of \textquotedblleft comma
category\textquotedblright\ is simply another way of describing bifunctors
$\mathcal{D}:\mathbf{X}^{op}\times\mathbf{A}\rightarrow\mathbf{Set}%
$.\footnote{Thus there are at least three ways to think of the elements of
$\mathcal{D}(x,a)$: as objects in the comma category $\mathcal{D}%
(1_{\mathbf{X}},1_{\mathbf{A}})$, as heteromorphisms $x\Rightarrow a$, or as
elements in a set-valued \textquotedblleft categorical
relation\textquotedblright\ from $\mathbf{X}$ to $\mathbf{A}$ \cite[p.
308]{bor:hca1} as in the profunctors-distributors-bimodules.} Starting with
the special case of $\mathcal{D}=Hom_{\mathbf{C}}(A(-),B(-))$, we get back the
comma category $(A,B)$ as $Hom_{\mathbf{C}}(A,B)$ (in the notation that makes
the bifunctor explicit). Starting with the hom-bifunctors $Hom_{\mathbf{A}%
}(F(-),1_{\mathbf{A}}(-))$ and $Hom_{\mathbf{X}}(1_{\mathbf{X}}(-),G(-))$,
this construction would reproduce Lawvere's comma categories $Hom_{\mathbf{A}%
}(F,1_{\mathbf{A}})=(F,1_{\mathbf{A}})$ and $Hom_{\mathbf{X}}(1_{\mathbf{X}%
},G)=(1_{\mathbf{X}},G)$.

Suppose that $Het(x,a)$ is represented on the left so that we have an
isomorphism $Hom_{\mathbf{A}}(Fx,a)\cong Het(x,a)$ natural in $x$ and in $a$.
Then the claim is that there is an isomorphism of comma categories
$Hom_{\mathbf{A}}(F,1_{\mathbf{A}})\cong Het(1_{\mathbf{X}},1_{\mathbf{A}})$
(using the notation with the bifunctor explicit) over the projections into the
product category. The representation associates a morphism $g(c):Fx\rightarrow
a$ with each element $x\overset{c}{\Rightarrow}a$ in $Het(x,a)$ so that it
provides the correspondence between the objects of the two comma categories.
Consider the two diagrams associated with morphisms in the two comma categories.

\begin{center}
$%
\begin{array}
[c]{ccccc}
& Fx & \overset{Fj}{\longrightarrow} & Fx^{\prime} & \\
g(c) & \downarrow &  & \downarrow & g(c^{\prime})\\
& a & \overset{k}{\longrightarrow} & a^{\prime} &
\end{array}
$ \ \ \ \ \ \ \ \ \ \ \ \ \ \ \ \ \ \ \ \ \ \ \ \ \ \ \ \ $%
\begin{array}
[c]{ccccc}
& x & \overset{j}{\longrightarrow} & x^{\prime} & \\
c & \Downarrow &  & \Downarrow & c^{\prime}\\
& a & \overset{k}{\longrightarrow} & a^{\prime} &
\end{array}
$
\end{center}

\noindent To see that the representation implies that one square will commute
if and only if the other does, consider the following commutative diagram
given by the representation.

\begin{center}
$%
\begin{array}
[c]{rcccl}%
g(c)\in & Hom_{\mathbf{A}}(Fx,a) & \cong & Het(x,a) & \backepsilon c\\
Hom(Fx,k) & \downarrow &  & \downarrow & Het(x,k)\\
& Hom_{\mathbf{A}}(Fx,a^{\prime}) & \cong & Het(x,a^{\prime}) & \\
Hom(Fj,a^{\prime}) & \uparrow &  & \uparrow & Het(j,a^{\prime})\\
g(c^{\prime})\in & Hom_{\mathbf{A}}(Fx^{\prime},a^{\prime}) & \cong &
Het(x^{\prime},a^{\prime}) & \backepsilon c^{\prime}%
\end{array}
$
\end{center}

\noindent Then $g(c)$ is carried by $Hom(Fx,k)$ to the same element in
$Hom_{\mathbf{A}}(Fx,a^{\prime})$ as $g(c^{\prime})$ is carried by
$Hom(Fj,a^{\prime})$ if and only if the square on the left above in the comma
category $Hom_{\mathbf{A}}(F,1_{\mathbf{A}})=(F,1_{\mathbf{A}})$ commutes.
Similarly, $c$ is carried by $Het(x,k)$ to the same element in
$Het(x,a^{\prime})$ as $c^{\prime}$ is carried by $Het(j,a^{\prime})$ if and
only if the square on the right above in the comma category $Het(1_{\mathbf{X}%
},1_{\mathbf{A}})$ commutes. Thus given the representation isomorphisms, the
square in the one comma category commutes if and only if the other one does.
With some more checking, it can be verified that the two comma categories are
then isomorphic over the projections.

\begin{center}
$%
\begin{array}
[c]{ccccc}
& Hom_{\mathbf{A}}(F,1_{\mathbf{A}}) & \cong & Het(1_{\mathbf{X}%
},1_{\mathbf{A}}) & \\
(\pi_{0},\pi_{1}) & \downarrow &  & \downarrow & (\pi_{0},\pi_{1})\\
& \mathbf{X}\times\mathbf{A} & = & \mathbf{X}\times\mathbf{A} &
\end{array}
$
\end{center}

\noindent Thus the `half-adjunction' in the language of representations is
equivalent to the `half-adjunction' in the language of comma categories. With
the other representation of $Het(x,a)$ on the right, we would then have the
comma category version of the three naturally isomorphic
bifunctors.\footnote{I owe to Colin McLarty the observation that the objects
of the comma categories $(F,1_{\mathbf{A}})$ and $(1_{\mathbf{X}},G)$ are
\textquotedblleft essentially the same\textquotedblright\ as the chimera
morphisms---as is shown by these isomorphisms.}

\begin{center}
$%
\begin{array}
[c]{ccccc}%
Hom_{\mathbf{A}}(F,1_{\mathbf{A}}) & \cong & Het(1_{\mathbf{X}},1_{\mathbf{A}%
}) & \cong & Hom_{\mathbf{X}}(1_{\mathbf{X}},G)\\
\downarrow &  & \downarrow &  & \downarrow\\
\mathbf{X}\times\mathbf{A} & = & \mathbf{X}\times\mathbf{A} & = &
\mathbf{X}\times\mathbf{A}%
\end{array}
$
\end{center}

\subsection{Adjunction Representation Theorem}

Adjunctions (particularly ones that are only abstractly specified) may not
have \textit{concretely} defined het-bifunctors that yield the adjunction via
the two representations. \ However, given any adjunction, there is always an
\textquotedblleft abstract\textquotedblright\ associated \ het-bifunctor given
by the main diagonal maps in the commutative adjunctive squares:
\[
Het(\widehat{x},\widehat{a})=\{\widehat{x}=(x,Fx)\overset{(f,f^{\ast}%
)}{\longrightarrow}(Ga,a)=\widehat{a}\}.
\]

\noindent The diagonal maps are closed under precomposition with maps from
$\widehat{\mathbf{X}}$ and postcomposition with maps from $\widehat
{\mathbf{A}}$. Associativity follows from the associativity in the ambient
category $\mathbf{X}\times\mathbf{A}$.

The representation is accomplished essentially by putting a $\widehat{hat}$ on
objects and morphisms embedded in $\mathbf{X}\times\mathbf{A}$. The categories
$\mathbf{X}$ and $\mathbf{A}$ are represented respectively by the subcategory
$\widehat{\mathbf{X}}$ with objects $\widehat{x}=(x,Fx)$ and morphisms
$\widehat{f}=(f,Ff)$ and by the subcategory $\widehat{\mathbf{A}}$ with
objects $\widehat{a}=(Ga,a)$ and morphisms $\widehat{g}=(Gg,g)$. The $(F,G)$
twist functor restricted to $\widehat{\mathbf{X}}\cong\mathbf{X}$ is
$\widehat{F}$ which has the action of $F$, i.e., $\widehat{F}\widehat
{x}=(F,G)(x,Fx)=(GFx,Fx)=\widehat{Fx}$ and similarly for morphisms. The twist
functor restricted to $\widehat{\mathbf{A}}\cong\mathbf{A}$ yields
$\widehat{G}$ which has the action of $G$, i.e., $\widehat{G}\widehat
{a}=(F,G)(Ga,a)=(Ga,FGa)=\widehat{Ga}$ and similarly for morphisms. These
functors provide representations on the left and right of the abstract
het-bifunctor $Het(\widehat{x},\widehat{a})=\{\widehat{x}\overset{(f,f^{\ast
})}{\longrightarrow}\widehat{a}\}$, i.e., the natural isomorphism%

\[
Hom_{\widehat{\mathbf{A}}}(\widehat{F}\widehat{x},\widehat{a})\cong
Het(\widehat{x},\widehat{a})\cong Hom_{\widehat{\mathbf{X}}}(\widehat
{x},\widehat{G}\widehat{a}).
\]

\ This birepresentation of the abstract het-bifunctor gives an isomorphic copy
of the original adjunction between the isomorphic copies $\widehat{\mathbf{X}%
}$ and $\widehat{\mathbf{A}}$ of the original categories. This theory of
adjoints may be summarized in the following:

\begin{description}
\item[\textbf{Adjunction Representation Theorem}:] Every adjunction
$F:\mathbf{X}\rightleftarrows\mathbf{A}:G$ can be represented (up to
isomorphism) as arising from the left and right representing universals of a
het-bifunctor $Het:\mathbf{X}^{op}\times\mathbf{A}\rightarrow\mathbf{Set}$
giving the chimera morphisms from the objects in a category $\mathbf{X}%
\cong\widehat{\mathbf{X}}$ to the objects in a category $\mathbf{A}%
\cong\widehat{\mathbf{A}}$.\footnote{In a historical note \cite[p.
103]{mac:cwm}, MacLane noted that Bourbaki \textquotedblleft
missed\textquotedblright\ the notion of an adjunction because Bourbaki focused
on the left representations of bifunctors $W:\mathbf{X}^{op}\times
\mathbf{A}\rightarrow\mathbf{Sets}$. MacLane remarks that given $G:\mathbf{A}%
\rightarrow\mathbf{X}$, they should have taken $W(x,a)=Hom_{\mathbf{X}}(x,Ga)$
and then focused on \textquotedblleft the symmetry of the adjunction
problem\textquotedblright\ to find $Fx$ so that $Hom_{\mathbf{A}}(Fx,a)\cong
Hom_{\mathbf{X}}(x,Ga)$. Thus MacLane missed the completely symmetrical
adjunction problem which is to find $Ga$ and $Fx$ such that $Hom_{\mathbf{A}%
}(Fx,a)\cong W(x,a)\cong Hom_{\mathbf{X}}(x,Ga)$.}
\end{description}

\subsection{Chimera Natural Transformations}

The concrete heteromorphisms, say, from sets to groups are just as `real' in
any relevant mathematical sense as the set-to-set functions or group-to-group
homomorphisms. At first, it might seem that the chimera morphisms do not fit
into a category theoretic framework since they do not compose with one another
like morphisms within a category. For instance, a set-to-functor cone cannot
be composed with another set-to-functor cone. But the proper sort of
composition was defined by having the morphisms of a category act on chimera
morphisms to yield other chimera morphisms. For instance, a set-to-set
function acts on a set-to-functor cone (as always, when composition is
defined) to yield another set-to-functor cone, and similarly for
functor-to-functor natural transformations composed on the other side. That
action on each side of a set of category-bridging heteromorphisms is described
by a het-bifunctor $Het:\mathbf{X}^{op}\times\mathbf{A}\rightarrow
\mathbf{Set}$. \ Thus the transition from the ordinary composition of
morphisms within a category to the sort of composition appropriate to
morphisms between categories is no more of a conceptual leap than going from
one bifunctor $Hom:\mathbf{X}^{op}\times\mathbf{X}\rightarrow\mathbf{Set}$ to
another $Het:\mathbf{X}^{op}\times\mathbf{A}\rightarrow\mathbf{Set}$. All
adjunctions arise as the birepresentations of such het-bifunctors.

Ordinary hom-bifunctors are (by definition) represented on the left and on the
right by the identity functor so the self-adjoint identity functor on any
category could be thought of as the \textquotedblleft
ur-adjunction\textquotedblright\ that expresses intra-category determination
in an adjunctive framework. Replace the birepresentable hom-bifunctor
$Hom:\mathbf{X}^{op}\times\mathbf{X}\rightarrow\mathbf{Set}$ with a
birepresentable het-bifunctor $Het:\mathbf{X}^{op}\times\mathbf{A}%
\rightarrow\mathbf{Set}$ and we have the smooth transition from intra-category
determination to cross-category determination via adjunctions.

Given the ubiquity of such adjoints in mathematics, there seems to be a good
case to stop treating chimera morphisms as some sort of `dark matter'
invisible to category theory; chimeras should be admitted into the `zoo'\ of
category theoretic creatures. Heteromorphisms (defined by the properties
necessary to be the elements of a het-bifunctor) between the objects of
different categories should be taken as entities in the ontology of category
theory at the same level as the morphisms between the objects within a
category (morphisms defined by the properties necessary to be the elements of
a hom-bifunctor).

New possibilities arise. For instance, the notion of a natural transformation
immediately generalizes to functors with different codomains by taking the
components to be heteromorphisms. \ Given functors $F:\mathbf{X}%
\rightarrow\mathbf{A}$ and $H:\mathbf{X}\rightarrow\mathbf{B}$ and a
het-bifunctor $Het:\mathbf{A}^{op}\times\mathbf{B}\rightarrow\mathbf{Set}$, a
\textit{chimera natural transformation }(\textit{relative to }$Het$),
$\varphi:F\Rightarrow H$, is given by set of morphisms $\{\varphi_{x}\in
Het(Fx,Hx)\}$ indexed by the objects of $\mathbf{X}$ such that for any
$j:x\rightarrow x^{\prime}$ the following diagram commutes

\begin{center}
$%
\begin{array}
[c]{ccccc}
& Fx & \overset{\varphi_{x}}{\Longrightarrow} & Hx & \\
Fj & \downarrow &  & \downarrow & Hj\\
& Fx^{\prime} & \overset{\varphi_{x^{\prime}}}{\Longrightarrow} & Hx^{\prime}
&
\end{array}
$
\end{center}

\noindent\lbrack the composition $Fx\overset{\varphi_{x}}{\Longrightarrow
}Hx\overset{Hj}{\longrightarrow}Hx^{\prime}$ is $Het(Fx,Hj)(\varphi_{x})\in
Het(Fx,Hx^{\prime})$, the composition $Fx\overset{Fj}{\longrightarrow
}Fx^{\prime}\overset{\varphi_{x^{\prime}}}{\Longrightarrow}Hx^{\prime}$ is
$Het(Fj,Hx^{\prime})(\varphi_{x^{\prime}})\in Het(Fx,Hx^{\prime})$, and
commutativity means they are the same element of $Het(Fx,Hx^{\prime})$]. These
chimera natural transformations do not compose like the morphisms in a functor
category but they are acted upon by the natural transformations in the functor
categories on each side to yield another bifunctor $Het^{X}%
(F,H)=\{F\Rightarrow H\}$.

There are chimera natural transformations each way between any functor and the
identity on its domain if the functor itself is used to define the appropriate
bifunctor $Het$. Given any functor $F:\mathbf{X}\rightarrow\mathbf{A}$, there
is a chimera natural transformation $1_{\mathbf{X}}\Rightarrow F$ relative to
the bifunctor defined as $Het(x,a)=Hom_{\mathbf{A}}(Fx,a)$ as well as a
chimera natural transformation $F\Rightarrow1_{\mathbf{X}}$ relative to
$Het(a,x)=Hom_{\mathbf{A}}(a,Fx)$.

Chimera natural transformations `in effect' already occur with reflective (or
coreflective) subcategories. A subcategory $\mathbf{A}$ of a category
$\mathbf{B}$ is a \textit{reflective subcategory} if the inclusion functor
$K:\mathbf{A}\rightarrow\mathbf{B}$ has a left adjoint. For any such
reflective adjunctions, the heteromorphisms $Het(b,a)$ are the $\mathbf{B}%
$-morphisms with their heads in the subcategory $\mathbf{A}$ so the
representation on the right $Het(b,a)\cong Hom_{\mathbf{B}}(b,Ka)$ is trivial.
The left adjoint $F:\mathbf{B}\rightarrow\mathbf{A}$ gives the representation
on the left: $Hom_{\mathbf{A}}(Fb,a)\cong Het(b,a)\cong Hom_{\mathbf{B}%
}(b,Ka)$. Then it is perfectly `natural' to see the unit of the adjunction as
defining a natural transformation $\eta:1_{\mathbf{B}}\Rightarrow F$ but that
is actually a chimera natural transformation (since the codomain of $F$ is
$\mathbf{A}$). Hence the conventional treatment (e.g., \cite[p. 89]{mac:cwm})
is to define another functor $R$ with the same domain and values on objects
and morphisms as $F$ except that its codomain is taken to be $\mathbf{B}$ so
that we can then `legally' have a natural transformation $\eta:1_{\mathbf{X}%
}\rightarrow R$ between two functors with the same codomain. Similar remarks
hold for the dual coreflective case where the inclusion functor has a right
adjoint and where the heteromorphisms are turned around, i.e., are
$\mathbf{B}$-morphisms with their tail in the subcategory $\mathbf{A}$.

The insertion of the generators maps $\{x\Rightarrow Fx\}$ define another
chimera natural transformation $h:1_{\mathbf{X}}\Rightarrow F$ from the
identity functor on the category of sets to the free-group functor (the
chimera version of the unit $\eta:1_{\mathbf{X}}\rightarrow GF$). A category
theory without chimeras can explain the naturality of maps such as
$x\rightarrow GFx$ but not the naturality of maps $x\Rightarrow Fx$!

For functors that are part of adjunctions, we can consider the embedding of
the adjunction representation theorem in the product category $\mathbf{X}%
\times\mathbf{A}$. Even if the adjunction is only abstractly given (so that we
have no concrete chimera maps), the representation theorem shows the maps
$(\eta_{x},1_{Fx}):\widehat{x}=(x,Fx)\rightarrow(GFx,Fx)=\widehat{Fx}$ (left
vertical arrows in adjunctive squares which are later shown to be
\textquotedblleft sending universals\textquotedblright) have the role of the
chimera units $h_{x}:x\Rightarrow Fx$. Moreover, they define a chimera natural
transformation $1_{\widehat{\mathbf{X}}}\Rightarrow\widehat{F}$ from the
identity functor on $\widehat{\mathbf{X}}$ to the functor $\widehat
{F}:\widehat{\mathbf{X}}\rightarrow\widehat{\mathbf{A}}$ (which is the $(F,G)
$ twist functor restricted to $\widehat{\mathbf{X}}$). Dually even if chimera
counits $e_{a}:Ga\Rightarrow a$ are not available, we still have the chimera
natural transformation $\widehat{G}\Rightarrow1_{\widehat{\mathbf{A}}}$ given
by the maps $(1_{Ga},\varepsilon_{a}):\widehat{G}\widehat{a}=\widehat
{Ga}=(Ga,FGa)\rightarrow(Ga,a)=\widehat{a}$ (right vertical arrows in
adjunctive squares which are later shown to be \textquotedblleft receiving
universals\textquotedblright) where $\widehat{G}:\widehat{\mathbf{A}%
}\rightarrow\widehat{\mathbf{X}}$ is the $(F,G)$ twist functor restricted to
$\widehat{\mathbf{A}}$. Thus every adjunction yields two chimera natural
transformations, the chimera versions of the unit and counit (we will see two
more later).

For the adjunctions of MacLane's `working mathematician,' there `should' be
concrete chimera morphisms that can be used to define these two chimera
natural transformations without resort to the representation theorem. In the
case of the limit adjunction, the chimera morphisms are \textquotedblleft
cones\textquotedblright\ from a set $w$ to a functor $D:\mathbf{D}%
\rightarrow\mathbf{Set}$ (where $\mathbf{D}$ is a diagram category). As with
every adjunction, there are two chimera natural transformations associated
with the unit and counit which in this case are $h:1_{\mathbf{Set}}%
\Rightarrow\Delta$ and $e:Lim\Rightarrow1_{\mathbf{Set}^{\mathbf{D}}}$ with
the chimera components given respectively by the two \textquotedblleft
universal cones\textquotedblright\ $w\Rightarrow\Delta w$ and $LimD\Rightarrow
D$ (see the section on limits in sets for more explanation).

While not `officially' acknowledged, the chimera morphisms or heteromorphisms
between objects of different categories are quite visible under ordinary
circumstances, once one acquires an eye to `see' them. This
non-acknowledgement is facilitated by the common practice of passing
effortlessly between the chimera morphism and its representation on one side
or the other. For instance in the last example, the chimera cone $w\Rightarrow
D$ from a set to a functor is often treated interchangeably with a natural
transformation $\Delta w\rightarrow D$ (its representation on the left) and
both are called \textquotedblleft cones.\textquotedblright\ We have argued
that chimera morphisms occur under quite ordinary circumstances, they fit
easily into a category-theoretic framework, and they spawn some new creatures
themselves---such as the chimera natural transformations.\footnote{This might
also have implications for n-category theory. The basic motivating example
from 1-category theory sees object-to-object morphisms only within a category
so that connections between categories are only by functors. \ And natural
transformations are seen only as going between two functors pinched together
at both ends to form a lens-shaped area. With chimera natural transformations,
the two functors only need to be pinched together on the domain side.}

\section{Chimera Factorizations}

\subsection{Abstract Adjunctive Squares}

These relationships can be conveniently and suggestively restated in the
adjunctive squares framework. \ Given $f:x\rightarrow Ga$ or $g:Fx\rightarrow
a$, the rest of the adjunctive square is determined so that it commutes and
$f=g^{\ast}$ and $g=f^{\ast}$. The two associated determinations are the two
components of the main diagonal $(f,g)$ so it can be thought of as one
abstract heteromorphism from $\widehat{x}=(x,Fx)$ to $\widehat{a}=(Ga,a)$.

\begin{center}
$%
\begin{array}
[c]{rcccl}
& (x,Fx) & \overset{(f,Ff)}{\rightarrow} & (Ga,FGa) & \\
(\eta_{x},1_{Fx}) & \downarrow & \nearrow & \downarrow & (1_{Ga}%
,\varepsilon_{a})\\
& (GFx,Fx) & \overset{(Gg,g)}{\longrightarrow} & (Ga,a) &
\end{array}
$
\end{center}

Of all the determinations $(x,Fx)\rightarrow(Ga,a)$ \textit{to} $(Ga,a)$, the
one that represents `self-determination' is the receiving universal
$(1_{Ga},\varepsilon_{a}):(Ga,FGa)\rightarrow(Ga,a)$ where $x=Ga$. All other
instances of a determination to $(Ga,a)$, e.g., $(GFx,Fx)\overset
{(Gg,g)}{\longrightarrow}(Ga,a)$, factor uniquely through the receiving
universal by the anti-diagonal map $(Gg,Fg^{\ast})$, i.e.,

\begin{center}
$(GFx,Fx)\overset{(Gg,g)}{\longrightarrow}(Ga,a)=(GFx,Fx)\overset
{(Gg,Fg^{\ast})}{\longrightarrow}(Ga,FGa)\overset{(1_{Ga},\varepsilon_{a}%
)}{\longrightarrow}(Ga,a)$.
\end{center}

Of all the determinations $(x,Fx)\rightarrow(Ga,a)$ \textit{from} $(x,Fx)$,
the one that represents `self-determination' is the sending universal
$(\eta_{x},1_{Fx}):(x,Fx)\rightarrow(GFx,Fx)$ where $a=Fx$. \ All other
instances of a determination from $(x,Fx)$, e.g., $(x,Fx)\overset
{(f,Ff)}{\longrightarrow}(Ga,FGa)$, factor uniquely through the sending
universal by the anti-diagonal map $(Gf^{\ast},Ff)$, i.e.,

\begin{center}
$(x,Fx)\overset{(f,Ff)}{\rightarrow}(Ga,FGa)=(x,Fx)\overset{(\eta_{x},1_{Fx}%
)}{\longrightarrow}(GFx,Fx)\overset{(Gf^{\ast},Ff)}{\longrightarrow}(Ga,FGa)$.
\end{center}

\noindent Since each of these factorizations of a top or bottom arrow goes
over to the other subcategory (e.g., $\widehat{\mathbf{X}}$ or $\widehat
{\mathbf{A}}$) and then back, they might be called the \textit{over-and-back
factorizations}.

Hence any $(x,Fx)$-to-$(Ga,a)$ heteromorphism given by the main diagonal
$(f,g)$ in a commutative adjunctive square (where $g=f^{\ast}$ and $f=g^{\ast
}$) factors through both the self-determination universals at the sending and
receiving ends by the anti-diagonal map $(Gg,Ff)$ obtained by applying the
$(F,G)$ twist functor to $\left(  f,g\right)  $. This factorization through
the two universals will be called the \textit{zig-zag factorization}. The
\textquotedblleft zig\textquotedblright\ ($\downarrow\nearrow$) followed by
the \textquotedblleft zag\textquotedblright\ ($\nearrow\downarrow$) give the
zig-zag ($\downarrow\nearrow\downarrow$) factorization.

\begin{description}
\item[Zig\noindent-Zag Factorization:] Any $(x,Fx)$-to-$(Ga,a)$ heteromorphism
in an adjunction factors uniquely through the sending and receiving universals
by the anti-diagonal map.
\end{description}

Recall that $G$ has to be one-to-one on morphisms of the form $g:Fx\rightarrow
a$ (uniqueness in the UMP for the counit) and $F$ has to be one-to-one on
morphisms of the form $f:x\rightarrow Ga$ (uniqueness in the UMP for the
unit). Thus the anti-diagonal maps of the form $(Gf^{\ast},Ff):\widehat
{Fx}\rightarrow\widehat{Ga}$ are uniquely correlated with $f$. Hence we can
define another bifunctor $Z(\widehat{F}\widehat{(-)},\widehat{G}\widehat
{(-)}):\mathbf{X}^{op}\times\mathbf{A}\rightarrow\mathbf{Set}$ of zig-zag
factorization maps:%

\[
\mathcal{Z}(\widehat{Fx},\widehat{Ga})=\{(Gf^{\ast},Ff):\widehat
{Fx}\rightarrow\widehat{Ga}\}.
\]

\noindent Since the anti-diagonal zig-zag factor map is uniquely determined by
the main diagonal map in a commutative adjunctive square, the $(F,G)$ twist
functor that takes $(f,f^{\ast})$ to $(Gf^{\ast},Ff)$ is an isomorphism from
$Het$ to $\mathcal{Z}$ which is easily checked to be natural in $x$ and in
$a$. Bearing in mind that these chimera bifunctors play the role of hom-sets
for chimera, this can be seen as another adjunction-like isomorphism between
the chimera morphisms in one direction and the chimera morphisms in the other direction:

\begin{center}
$Het(\widehat{x},\widehat{a})\cong\mathcal{Z}(\widehat{Fx},\widehat{Ga})$.
\end{center}

\noindent In spite of the notation $\mathcal{Z}(\widehat{Fx},\widehat{Ga})$,
there is no implication that we have a bifunctor $\mathcal{Z}(\widehat
{a},\widehat{x})$\ defined on arbitrary $\widehat{a}$ and $\widehat{x}$, and
thus no implication that we are dealing with arbitrary chimera morphisms
$\widehat{Fx}\Rightarrow\widehat{Ga}$. The bifunctor $\mathcal{Z}(\widehat
{Fx},\widehat{Ga})$ is here \textit{defined} as the image of $Het(\widehat
{x},\widehat{a})$ under the $(F,G)$ twist functor.

In the representation theorem, we saw that every adjunction has a
het-bifunctor $Het(\widehat{x},\widehat{a})$ of chimera morphisms from
$\widehat{\mathbf{X}}$ to $\widehat{\mathbf{A}}$ so that the adjunction arises
(up to isomorphism) as the birepresentation of that het-bifunctor
$Hom_{\widehat{\mathbf{A}}}(\widehat{F}\widehat{x},\widehat{a})\cong
Het(\widehat{x},\widehat{a})\cong Hom_{\widehat{\mathbf{X}}}(\widehat
{x},\widehat{G}\widehat{a})$. Now we see that every adjunction also gives rise
to another bifunctor $\mathcal{Z}(\widehat{F}\widehat{x},\widehat{G}%
\widehat{a})$ of heteromorphisms going in the other direction from
$\widehat{\mathbf{A}}$ to $\widehat{\mathbf{X}}$ and that the two bifunctors
are naturally isomorphic. Thus we now have four different bifunctors that are
naturally isomorphic in $x$ and in $a$.

\begin{center}
$Hom_{\widehat{\mathbf{A}}}(\widehat{F}\widehat{x},\widehat{a})\cong
Het(\widehat{x},\widehat{a})\cong\mathcal{Z}(\widehat{F}\widehat{x}%
,\widehat{G}\widehat{a})\cong Hom_{\widehat{\mathbf{X}}}(\widehat{x}%
,\widehat{G}\widehat{a})$
\end{center}

\noindent Given $f:a\rightarrow Ga$ and thus $f^{\ast}:Fx\rightarrow a$, these
isomorphisms correlate together four morphism which in the adjunctive square
are the bottom map, the main diagonal, the anti-diagonal, and the top map.

In an adjunction, the operation of taking the adjoint transpose amounts to
applying the appropriate functor and pre- or post-composing with the
appropriate universal. For $f:x\rightarrow Ga$, $f^{\ast}=\varepsilon_{a}Ff $
and for $g:Fx\rightarrow a$, $g^{\ast}=Gg\eta_{x}$. For the chimera
isomorphism, the transpose of $\widehat{x}=(x,Fx)\overset{(f,f^{\ast}%
)}{\longrightarrow}(Ga,a)=\widehat{a}$ is obtained simply by applying the
$(F,G)$ twist functor to get $\widehat{Fx}=(GFx,Fx)\overset{(Gf^{\ast}%
,Ff)}{\longrightarrow}(Ga,FGa)=\widehat{Ga}$. Taking the transpose in the
other direction is done by pre- and post-composing both the universals, i.e.,
the zig-zag factorization. Both isomorphisms are illustrated in a commutative
adjunctive square. The adjoint transposes are the top and bottom arrows and
the chimera adjoint correlates are the main diagonal and the anti-diagonal.

\subsection{Abstract Adjunctive-Image Squares}

It was previously noted that the uniqueness requirement in the UMPs of an
adjunction imply that the functor $F$ has to be one-one on morphisms of the
form $f:x\rightarrow Ga$ while $G$ has to be one-one on morphisms
$g:Fx\rightarrow a$. Thus the $(F,G)$ twist functor has to be one-one on all
the morphisms in an adjunctive square. Applying the twist functor yields the
image of the (commutative) adjunctive square which we will call the
\textit{adjunctive-image square}.

\begin{center}
$%
\begin{array}
[c]{rcccl}%
\widehat{Fx}= & (GFx,Fx) & \overset{(GFf,Ff)}{\longrightarrow} & (GFGa,FGa) &
=\widehat{FGa}\\
(1_{GFx},F\eta_{x}) & \downarrow &  & \downarrow & (G\varepsilon_{a}%
,1_{FGa})\\
\widehat{GFx}= & (GFx,FGFx) & \overset{(Gf^{\ast},FGf^{\ast})}{\longrightarrow
} & (Ga,FGa) & =\widehat{Ga}%
\end{array}
\medskip$

Abstract Adjunctive-Image Square
\end{center}

\noindent The top map is in $\widehat{\mathbf{A}}$, the bottom map is in
$\widehat{\mathbf{X}}$, and the main diagonal is the anti-diagonal of the
adjunctive square. Again there is a unique anti-diagonal map (obtained again
as the image under the $(F,g)$ twist functor) in the adjunctive-image square,
i.e., $(GFf,FGf^{\ast}):(GFx,FGFx)\rightarrow(GFGa,FGa)$, which makes the
upper and lower triangles commute. Hence the original anti-diagonal map
$(Gf^{\ast},Ff)$ has a unique zig-zag factorization through this new
anti-diagonal map using the two vertical universal maps.

The image of the hom-sets such as $Hom_{\widehat{\mathbf{X}}}(\widehat
{x},\widehat{G}\widehat{a})$ under $\widehat{F}$ will yield an isomorphic set
of morphisms like the top maps in the adjunctive-image square but it is not
necessarily a hom-set. This isomorphic image is a bifunctor which we could
denote $H(\widehat{F}\widehat{x},\widehat{FG}\widehat{a})$ where there is a
natural isomorphism $Hom_{\widehat{\mathbf{X}}}(\widehat{x},\widehat
{G}\widehat{a})\cong H(\widehat{F}\widehat{x},\widehat{FG}\widehat{a})$.
\ Similarly, we have $Hom_{\widehat{\mathbf{A}}}(\widehat{F}\widehat
{x},\widehat{a})\cong H(\widehat{GF}\widehat{x},\widehat{G}\widehat{a})$
concerning the bottom maps in the adjunctive-image square. This use of the
$(F,G)$ twist functor to obtain an isomorphic bifunctor from a given one was
already used to derive the bifunctor of anti-diagonal maps from the bifunctor
of diagonal maps: $Het(\widehat{x},\widehat{a})\cong\mathcal{Z}(\widehat
{F}\widehat{x},\widehat{G}\widehat{a})$.\footnote{The anti-diagonal in the
adjunctive square becomes the main diagonal in the adjunctive-image square.
That process can be repeated so that \textquotedblleft it's chimeras all the
way down.\textquotedblright} In all three cases, the process can be continued
to derive an infinite sequence of bifunctors all isomorphic to one another.
Only the first stage in the sequence for the chimera maps $Het(\widehat
{x},\widehat{a})\cong\mathcal{Z}(\widehat{F}\widehat{x},\widehat{G}\widehat
{a})$ will be investigated here.

\subsection{Chimera Adjunctive Squares}

In any concretely specified adjunction (i.e., not just the abstract definition
used in the representation theorem) that occurs in mathematics, we would
\textit{expect} to be able to `takes the hats off' and find concrete chimera
bifunctors $Het(x,a)$ and $\mathcal{Z}(Fx,Ga)$ to give us the same natural isomorphisms:%

\[
Hom_{\mathbf{A}}(Fx,a)\cong Het(x,a)\cong\mathcal{Z}(Fx,Ga)\cong
Hom_{\mathbf{X}}(x,Ga).
\]

\noindent In this section, we give the relationships that the abstract
adjunctive square tells us to expect to find for the heteromorphisms in the
two chimera bifunctors. In each example, we need to find these bifunctors and
demonstrate these relationships.

Suppose we have the het-bifunctor, the `zig-zag' bifunctor and the
birepresentations to give the above isomorphisms. We previously used the
representations of $Het(x,a)$ to pick out universal elements, the chimera unit
$h_{x}\in Het(x,Fx)$ and the chimera counit $e_{a}\in Het(Ga,a)$, as the
respective correlates of $1_{Fx}$ and $1_{Ga}$ under the isomorphisms. In the
isomorphisms of the four bifunctors, let $g(c):Fx\rightarrow a$,
$c:x\Rightarrow a$, $z(c):Fx\Rightarrow Ga$, and $f(c):x\rightarrow Ga$ be the
four adjoint correlates. \ We showed that from the birepresentation of
$Het(x,a)$, any chimera morphism $x\overset{c}{\Rightarrow}a$ in $Het(x,a)$
would have two factorizations: $g(c)h_{x}=c=e_{a}f(c)$. This two
factorizations are spliced together along the main diagonal $c:x\Rightarrow a
$ to form the chimera (commutative) adjunctive square.

\begin{center}
$%
\begin{array}
[c]{ccccc}
& x & \overset{f(c)}{\longrightarrow} & Ga & \\
h_{x} & \Downarrow &  & \Downarrow & e_{a}\\
& Fx & \overset{g(c)}{\longrightarrow} & a &
\end{array}
\medskip$

Chimera Adjunctive Square
\end{center}

\noindent In the examples, we need to find the chimera morphisms
$Fx\overset{z(c)}{\Longrightarrow}Ga$ that give us $\mathcal{Z}(Fx,Ga)$ and
that fit as the anti-diagonal morphisms in these chimera adjunctive squares to
give commutative upper and lower triangles. The two commutative triangles
formed by the anti-diagonal might be called the:

\begin{center}
$x\overset{h_{x}}{\Longrightarrow}Fx\overset{z(c)}{\Rightarrow}Ga=x\overset
{f(c)}{\longrightarrow}Ga$

$Fx\overset{z(c)}{\Rightarrow}Ga\overset{e_{a}}{\Longrightarrow}%
a=Fx\overset{g(c)}{\longrightarrow}a\medskip$

\textit{Over-and-back factorizations of }$f$\textit{\ and }$g$.
\end{center}

\noindent One of the themes of this theory of adjoints is that some of the
rigmarole of the conventional treatment of adjoints (\textit{sans} chimeras)
is only necessary because of the restriction to morphisms within one category
or the other. The chimera unit $x\overset{h_{x}}{\Longrightarrow}Fx $ only
involves one of the functors and it appears already with the half-adjunction
of a representation on the left, and similarly for the chimera counit
$Ga\overset{e_{a}}{\Longrightarrow}a$. The importance of the unit and counit
lies in their universal mapping properties for morphisms of the form
$x\overset{f}{\longrightarrow}Ga$ or $Fx\overset{g}{\longrightarrow}a.$ With
the over-and-back factorizations, we again see the power of adjoints expressed
in simpler terms using the underlying heteromorphisms. The same holds for the
triangular identities as we will see in the next section.

\subsection{Chimera Adjunctive-Image Squares}

When we have the four isomorphic bifunctors, then $1_{Fx}$ would associate
with a chimera denoted $h_{x2}=z(h_{x}):Fx\Rightarrow GFx$ in $\mathcal{Z}%
(Fx,GFx)$ and $1_{Ga}$ would associate with $e_{a1}=z(e_{a}):FGa\Rightarrow
Ga$ in $\mathcal{Z}(FGa,Ga)$. The notation is chosen since taking
$f(c)=\eta_{x}$ (namely $c=h_{x}$), we have the over-and-back factorization of
the unit: $x\overset{h_{x}}{\Longrightarrow}Fx\overset{h_{x2}}{\Longrightarrow
}GFx=x\overset{\eta_{x}}{\longrightarrow}GFx$, so $h_{x2}$ is the `second
part' post-composed to $h_{x}$ to yield $\eta_{x}$. And taking
$g(c)=\varepsilon_{a}$ (namely $c=e_{a}$), we have the over-and-back
factorization of the counit: $FGa\overset{e_{a1}}{\Longrightarrow}%
Ga\overset{e_{a}}{\Longrightarrow}a=FGa\overset{\varepsilon_{a}}%
{\longrightarrow}a$, so $e_{a1}$ is the `first part' pre-composed to $e_{a}$
to yield $\varepsilon_{a}$. \ 

These special chimera morphisms have universality properties since they are
associated with identity maps $1_{Fx}$ and $1_{Ga}$ in the above natural
isomorphism of four bifunctors (mimic the proofs of the universality
properties of $\eta_{x}$, $\varepsilon_{a}$, $h_{x}$, and $e_{a}$). The
`second half of the unit' $h_{x2}=z(h_{x})$ has the following universality
property: for any anti-diagonal map from $Fx$, $z(c):Fx\Rightarrow Ga$, there
is unique map $g=g(c):Fx\rightarrow a$ such that%

\[
Fx\overset{z(c)}{\Longrightarrow}Ga=Fx\overset{h_{x2}}{\Longrightarrow
}GFx\overset{Gg}{\longrightarrow}Ga.
\]

\noindent Precompose $h_{x}:x\Rightarrow Fx$ on both sides and we have the
usual UMP for the unit, i.e.,%

\[
x\overset{f(c)}{\longrightarrow}Ga=x\overset{h_{x}}{\Rightarrow}%
Fx\overset{z(c)}{\Rightarrow}Ga=x\overset{h_{x}}{\Rightarrow}Fx\overset
{h_{x2}}{\Rightarrow}GFx\overset{Gg(c)}{\longrightarrow}Ga=x\overset{\eta_{x}%
}{\longrightarrow}GFx\overset{Gf^{\ast}}{\longrightarrow}Ga.
\]

The same chimera morphism $z(c):Fx\Rightarrow Ga$ also has a factorization
through the other anti-diagonal universal, the first half of the counit,
$e_{a1}=z(e_{a})$. Given any anti-diagonal map to $Ga$, $z(c):Fx\Rightarrow Ga
$, there is a unique map $f=f(c):x\rightarrow Ga$ such that%

\[
Fx\overset{z(c)}{\Longrightarrow}Ga=Fx\overset{Ff}{\Longrightarrow}%
FGa\overset{e_{a1}}{\Longrightarrow}Ga.
\]

\noindent Post-composing with $e_{a}:Ga\Rightarrow a$ on both sides yields the
usual UMP for the counit, i.e.,%

\[
Fx\overset{g(c)}{\longrightarrow}a=Fx\overset{z(c)}{\Longrightarrow}%
Ga\overset{e_{a}}{\Rightarrow}a=Fx\overset{Ff}{\Longrightarrow}FGa\overset
{e_{a1}}{\Longrightarrow}Ga\overset{e_{a}}{\Rightarrow}a=Fx\overset{Fg^{\ast}%
}{\Longrightarrow}FGa\overset{\varepsilon_{a}}{\longrightarrow}a.
\]

Splicing together the two factorizations of $z(c):Fx\Rightarrow Ga$ as the
main diagonal, we have the chimera version of the adjunctive-image square.

\begin{center}
$%
\begin{array}
[c]{ccccc}
& Fx & \overset{Ff}{\longrightarrow} & FGa & \\
h_{x2} & \Downarrow &  & \Downarrow & e_{a1}\\
& GFx & \overset{Gf^{\ast}}{\longrightarrow} & Ga &
\end{array}
\medskip$

Chimera Adjunctive-Image Square
\end{center}

\noindent Thus the two universal properties give the two ways the main
diagonal $Fx\overset{z(c)}{\Longrightarrow}Ga$ factors through the two
vertical universal arrows $Fx\overset{h_{x2}}{\Longrightarrow}GFx$ and
$FGa\overset{e_{a1}}{\Longrightarrow}Ga$. As before with $1_{\mathbf{X}%
}\overset{h}{\Longrightarrow}F$ and $G\overset{e}{\Longrightarrow
}1_{\mathbf{A}}$, the chimera universals are also the components of two
chimera natural transformations: $F\overset{h_{2}}{\Rightarrow}GF$ and
$FG\overset{e_{1}}{\Rightarrow}G$. Thus every adjunction has associated with
it four chimera natural transformations, and the two conventional natural
transformations associated with an adjunction are obtained as composites of
the four chimera natural transformations.

\begin{center}
$1_{\mathbf{X}}\overset{h}{\Longrightarrow}F\overset{h_{2}}{\Longrightarrow
}GF=1_{\mathbf{X}}\overset{\eta}{\longrightarrow}GF$ and $FG\overset{e_{1}%
}{\Longrightarrow}G\overset{e}{\Longrightarrow}1_{\mathbf{A}}=FG\overset
{\varepsilon}{\longrightarrow}1_{\mathbf{A}}$.
\end{center}

Since the anti-diagonal $z(c)$ can be factored, the over-and-back
factorizations for $f:x\rightarrow Ga$ and for $g:Fx\rightarrow a$ can be
factored again as can be pictured using the adjunctive-image square. Add
$x\overset{h_{x}}{\Longrightarrow}Fx$ on as a pendant to the above chimera
adjunctive-image square (where $z(c)$ is the main diagonal) to obtain the
following diagram.

\begin{center}
$%
\begin{array}
[c]{cccccc}%
x & \overset{h_{x}}{\Longrightarrow} & Fx & \overset{Ff}{\longrightarrow} &
FGa & \\
& h_{x2} & \Downarrow &  & \Downarrow & e_{a1}\\
&  & GFx & \overset{Gf^{\ast}}{\longrightarrow} & Ga &
\end{array}
$
\end{center}

\noindent The pendant followed by the counter-clockwise maps gives the usual
factorization $x\overset{\eta_{x}}{\longrightarrow}GFx\overset{Gf^{\ast}%
}{\longrightarrow}Ga=x\overset{f}{\longrightarrow}Ga$. Hence following the
pendant by the clockwise maps gives a different \textit{over-across-and-back
factorization of }$f$. Dually, we could postcompose the pendant $Ga\overset
{e_{a}}{\Longrightarrow}a$ on the bottom and obtain the over-across-and-back
factorization of $g=f^{\ast}:Fx\rightarrow a$ in addition to the usual one.
These four factorizations may be summarized as follows:%

\begin{align*}
x\overset{h_{x}}{\Longrightarrow}Fx\overset{h_{x2}}{\Longrightarrow
}GFx\overset{Gf^{\ast}}{\longrightarrow}Ga  &  =x\overset{f}{\longrightarrow
}Ga\\
x\overset{h_{x}}{\Longrightarrow}Fx\overset{Ff}{\longrightarrow}%
FGa\overset{e_{a1}}{\Rightarrow}Ga  &  =x\overset{f}{\longrightarrow}Ga\\
Fx\overset{Fg^{\ast}}{\longrightarrow}FGa\overset{e_{a1}}{\Rightarrow
}Ga\overset{e_{a}}{\Longrightarrow}a  &  =Fx\overset{g}{\longrightarrow}a\\
Fx\overset{h_{x2}}{\Longrightarrow}GFx\overset{Gg}{\longrightarrow}%
Ga\overset{e_{a}}{\Longrightarrow}a  &  =Fx\overset{g}{\longrightarrow}a.
\end{align*}

Specializing $f=1_{Ga}$ gives one triangular identity in the first equation.
But in the second equation, it gives the first \textit{over-and-back identity}
(a `short form' of the triangular identity):%

\[
Ga\overset{h_{Ga}}{\Longrightarrow}FGa\overset{e_{a1}}{\Longrightarrow
}Ga=Ga\overset{1_{Ga}}{\longrightarrow}Ga.
\]

\noindent Specializing $g=1_{Fx}$ gives the other triangular identity in the
third equation. But in the fourth equation, it gives the other over-and-back identity:%

\[
Fx\overset{h_{x2}}{\Longrightarrow}GFx\overset{e_{Fx}}{\Longrightarrow
}Fx=Fx\overset{1_{Fx}}{\longrightarrow}Fx.
\]

\noindent Here again we see chimera natural transformations composing to yield
conventional natural transformations:

\begin{center}
$G\overset{h_{G}}{\Longrightarrow}FG\overset{e_{1}}{\Longrightarrow
}G=G\overset{1_{G}}{\longrightarrow}G$

$F\overset{h_{2}}{\Longrightarrow}GF\overset{e_{F}}{\Longrightarrow
}F=F\overset{1_{F}}{\longrightarrow}F.\medskip$

\textit{Over-and-back identities}
\end{center}

\noindent Thus on the functorial images $Ga$ and $Fx$, there is a canonical
heteromorphism to the other category and a canonical heteromorphism coming
back so that the composition is the identity on the images.

\section{Limits in Sets}

Category theory is about the horizontal transmission of structure or, more
generally, determination between objects. The important examples in
mathematics have the structure of determination through (self-participating)
universals. Such a universal has the property in question and then every other
instance of the property is determined to have it by \textquotedblleft
participating in\textquotedblright\ (e.g., uniquely factoring through) the
self-participating universal. Adjunction extends this theme so that the
determination has symmetrically both a sending universal and a receiving
universal and where all determinations factor through both universals by the
anti-diagonal morphism in an adjunctive square. The `self-determination'
involved in the universals was previously illustrated using the example of the
product and coproduct construction in sets. \ We now generalize the
illustration by considering the adjunction for limits in the category of sets
(and colimits in the next section).

The construction of limits in sets generalizes the example of the product. Let
$\mathbf{D}$ be a small (diagram) category and $D:\mathbf{D}{\rightarrow} $
$\mathbf{Set}$ a functor considered as a diagram in the category of
$\mathbf{Set}$. The diagram $D$ is in the functor category $\mathbf{Set}%
^{\mathbf{D}}$ where the morphisms are natural transformation between the
functors. Let $\Delta:\mathbf{Set}{\rightarrow}\mathbf{Set}^{\mathbf{D}}$, the
diagonal functor, assign to each set $w$ the constant functor $\Delta w$ on
$\mathbf{D}$ whose value for each $i$ in $\mathbf{D}$ is $w$, and for each
morphism $\alpha:i\rightarrow j$ in $\mathbf{D}$, the value of $\Delta
w_{\alpha}$ is $1_{w}$. The functor in the other direction is the limit
functor $Lim$ which assigns a set $LimD$ to each diagram $D$ and a set
function $Lim\theta:LimD\rightarrow LimD^{\prime}$ to every natural
transformation $\theta:D\rightarrow D^{\prime}$.

The adjunction for the diagonal and limit functors is:

\begin{center}
${Hom(\Delta w,D)\cong Hom(w,LimD)}$
\end{center}

\noindent where is $\Delta$ is the left adjoint and $Lim$ is the right
adjoint. A adjunctive square for this adjunction would have the generic form:

\begin{center}
$%
\begin{array}
[c]{ccc}%
(w,\Delta w) & \overset{(f,\Delta f)}{\longrightarrow} & (LimD,\Delta LimD)\\
\downarrow & \overset{^{(Limg,\Delta f)}}{\nearrow} & \downarrow\\
(Lim\Delta w,\Delta w) & \overset{(Limg,g)}{\longrightarrow} & (LimD,D)
\end{array}
\medskip$

Abstract Adjunctive Square for Limits Adjunction
\end{center}

\noindent which commutes when $g=f^{\ast}$or $f=g^{\ast}$.

For this adjunction, a chimera morphism or heteromorphism from a set $w$ to a
diagram functor $D$ is concretely given by a \textit{cone} $w\Rightarrow D$
which is defined as a set of maps $\{w\overset{f_{i}}{\longrightarrow}D_{i}\}$
indexed by the objects $i$ in the diagram category $\mathbf{D}$ such that for
any morphism $\alpha:i\rightarrow j$ in $\mathbf{D}$, $w\overset{f_{i}%
}{\longrightarrow}D_{i}\overset{D_{\alpha}}{\longrightarrow}D_{j}%
=w\overset{f_{j}}{\longrightarrow}D_{j}$. The adjunction is then given by the
birepresentation of the het-bifunctor where $Het(w,D)=\{w\Rightarrow D\}$ is
the set of cones from the set $w$ to the diagram functor $D$. Instead of
proceeding formally from the het-bifunctor, we will conceptually analyze the
construction of the functors $G=Lim$ and $F=\Delta$.

Conceptually, start with the idea of a set function as a way for elements in
the domain to determine certain `elements' in the codomain, the root concept
abstracted to form the concept of a morphism in category theory. In the case
of a diagram functor $D$ as the target for the determination, the conceptual
atom or \textit{element\ of }$D$ is the maximal set of points that could be
determined in the codomains $D_{i}$ by functions from a one point set and that
are compatible with the morphisms between the $D_{i}$s. Thus an
\textquotedblleft element\textquotedblright\ of $D$ is classically a
\textquotedblleft1-cone\textquotedblright\ which consists of a member $x_{i} $
of the set $D_{i}$ for each object $i$ in the diagram category $\mathbf{D}$
such that for every morphism $\alpha:i\rightarrow j$ in $\mathbf{D}$,
$D_{\alpha}(x_{i})=x_{j}$. Ordinarily, a element of $D$ would be thought of as
an element in the product of the $D_{i}$s that is compatible with the
morphisms, and $LimD$ would be the subset of compatible elements of the direct
product of all the $D_{i}$s.

An element of $D$ is the abstract `determinee' of a point \textit{before} any
point is assigned to be the `determiner'. Thus the notion of a 1-cone\ isn't
quite right conceptually in that it pictures the element as having been
already determined by a one point set. The notion of a global
section\ $\cite[p. 47]{macm:sh}$ is a better description that is free of this
connotation. The right adjoint represents these atomic determinees as an
object $LimD$ in the top category, in this case $\mathbf{Set}$. Then a
morphism in that category $f:w\rightarrow LimD$ will map one-point determiners
to the atomic determinees in $LimD$ and will thus give a determination from
the set $w$ to the diagram functor $D$ (i.e., a cone).

It is a conceptual move to take a `determinee' as the abstract `determiner' of
itself, and that `self-determination' by the map $1_{LimD}:LimD\rightarrow
LimD$ takes the set of determinees of $D$ reconceptualized as `determiner'
elements of a set. But to be a determiner by a morphism in the bottom
category, one needs a morphism in that functor category, namely a natural
transformation. A natural transformation $\theta:D^{\prime}\rightarrow D$
carries the determinees or elements of a functor $D^{\prime} $ to the elements
of the functor $D$. The determinees in the set $LimD$ are repackaged in the
diagram $\Delta LimD$ (see below) as determiners and the natural
transformation by which they canonically determine the determinees of the
diagram $D$ (from which they came) is the counit $\varepsilon_{D}:\Delta
LimD\rightarrow D$.

The conventional treatment of an adjunction is complicated by the need to deal
only with morphisms within either of the categories. That need also accounts
for the intertwining of the two functors in the unit and counit universals.
The properties of the chimera versions of the unit and counit follow from each
of the half-adjunctions (the representations on the left or right of the
het-bifunctor) without any intertwining of the two functors. For instance,
there is no composition of functors in the chimera adjunctive square (see
below). The additional relationships obtained from both representations and
the intertwining of the functors are illustrated in the adjunctive-image
square. For instance, the chimera version of the counit $\varepsilon_{D}$ is
$e_{D}:LimD\Rightarrow D$ is simply the cone of projection maps and it does
not involve the diagonal functor $\Delta$. Moreover, the role of
self-determination is much clearer. Thinking of the elements of $D$ as the
determinees, they become the determiners as the elements of $LimD$ and that
self-determination is given by the projection maps $e_{D}:LimD\Rightarrow D$.

We need to do a similar conceptual analysis of the other functor $\Delta$. How
can the elements of a set $w$ be determiners of some diagram functor on the
diagram category $\mathbf{D}$? The constant functor $\Delta w$ repackages the
set $w$ as the constant value of a functor on the diagram category. Since the
maps between the values of $\Delta w$ are all identity maps $1_{w}$, the
elements on each connected component of the diagram category are just the
elements of $w$. If the diagram category is multiply connected, then the
elements of any diagram $D$ on that category will be the product of the
elements of the functor restricted to each connected component. If the diagram
category $\mathbf{D}$ had, say, two components
$\mathbf{D\mathbf{\upharpoonright}1}$ and $\mathbf{D\mathbf{\upharpoonright}2}
$, then the elements of the functor $D$ would be ordered pairs of the elements
of $D\mathbf{\upharpoonright}1$ and $D\mathbf{\upharpoonright}2$, the functor
$D$ restricted to each component of the diagram category. Thus the set $LimD$
would have the structure of a product of sets of elements of the functors
$D\mathbf{\upharpoonright}1$ and $D\mathbf{\upharpoonright}2$, i.e.,
$LimD=LimD\mathbf{\upharpoonright}1\times LimD\mathbf{\upharpoonright}2 $.

A determination by a natural transformation $g:\Delta w\rightarrow D$ could be
parsed as two natural transformations $g_{k}:\Delta w\mathbf{\upharpoonright
}k\rightarrow D\mathbf{\upharpoonright}k$ for $k=1,2$, each one carrying an
element of $w$ to a determinee or element of $D\mathbf{\upharpoonright}k$.
Thus a determination or cone $w\Rightarrow D$ represented in the form
$g:\Delta w\rightarrow D$ would induce a determination $g^{\ast}:w\rightarrow
LimD=LimD\mathbf{\upharpoonright}1\times LimD\mathbf{\upharpoonright}2$ in the
other form. The constant functor $\Delta$ repackages a set $w$ as a diagram
functor to be a determiner of diagram functors. \ The set $w$ can also be
self-determining by taking $g=1_{\Delta w}:\Delta w\rightarrow\Delta w$. \ But
for that determination to be given by a morphism in the top category, the
determinees of $\Delta w$ are represented in the top category by applying the
limit functor. \ There is then a canonical (diagonal) map wherein $w$
determines the determinees $Lim\Delta w$ of its own representation as a
determiner $\Delta w$, namely $\eta_{w}:w\rightarrow Lim\Delta w$ (where with
two connected components $Lim\Delta w\cong w\times w$). As will be seen, the
chimera version is $h_{w}:w\Rightarrow\Delta w$ is simply the cone of identity
maps $1_{w}$.

In the adjunctions encountered by the working mathematician, there should be
concrete chimeras or heteromorphisms. This is true for all the examples
considered here. As noted above, the heteromorphisms from a set to a diagram
functor are the cones $w\Rightarrow D$ (\textit{Nota bene}, not the natural
transformations $\Delta w\rightarrow D$ which represent-on-the-left the
chimera cones $w\Rightarrow D$). When the concrete chimeras are available,
then the adjunctive square can be developed using them (instead of only the
abstract version embedded in the product category). The chimera version
$h_{w}:w\Rightarrow\Delta w$ of the unit is the cone where each function
$w\rightarrow(\Delta w)_{i}=w$ is $1_{w}$. The chimera version $e_{D}%
:LimD\Rightarrow D$ of the counit is the cone of projection maps. Given any
cone $c:w\Rightarrow D$, there is a unique set map $f(c):w\rightarrow LimD$
such that $w\overset{f(c)}{\longrightarrow}LimD\overset{e_{D}}{\Rightarrow
}D=w\overset{c}{\Rightarrow}D$.\footnote{We leave all the routine checking of
the chimera relationships to the reader.} And there is a unique natural
transformation $g(c):\Delta w\rightarrow D$ such that $w\overset{h_{w}%
}{\Rightarrow}\Delta w\overset{g(c)}{\longrightarrow}D=w\overset
{c}{\Rightarrow}D$. These mappings provide the two representations:%

\[
Hom(\Delta w,D)\cong Het(w,D)\cong Hom(w,LimD)
\]
These maps also give the chimera version of the adjunctive square where
$c:w\Rightarrow D$ is the main diagonal.

\begin{center}
$%
\begin{array}
[c]{ccccc}
& w & \overset{f(c)}{\longrightarrow} & LimD & \\
h_{w} & \Downarrow &  & \Downarrow & e_{D}\\
& \Delta w & \overset{g(c)}{\longrightarrow} & D &
\end{array}
\medskip$

Chimera Adjunctive Square for Limits Adjunction
\end{center}

\noindent The anti-diagonal map for the zig-zag factorization will be a
chimera morphism going from the diagram $\Delta w$ to the set $LimD$ (see
cocones in the next section). For each $i$ in the diagram category, $(\Delta
w)_{i}=w$ so each component of the \textquotedblleft cocone\textquotedblright%
\ $z(c):\Delta w\Rightarrow LimD$ is $f(c):w\rightarrow LimD$. Often the
anti-diagonal chimera can be defined using a heteromorphic inverse to one of
the vertical morphisms. In this case, $h_{w}:w\Rightarrow\Delta w$ has\ the
inverse cocone $\Delta w\Rightarrow w$ which has the identity map $1_{w}$ as
each component. Joining the cone and cocone at the open end just yields the
identity: $w\overset{h_{w}}{\Rightarrow}\Delta w\Rightarrow w=1_{w}$ and
joining them at the vertexes yields the other identity $\Delta w\Rightarrow
w\overset{h_{w}}{\Rightarrow}\Delta w=1_{\Delta w}$.\footnote{Composing a cone
and cocone (same diagram category) on the open end to yield a set map works
intuitively if the diagram category is connected or if all the maps in the
cone are the same and similarly in the cocone (the latter being the case
here). Note that we have here the heteromorphic version of an `isomorphism'
between a set $w$ and a functor $\Delta w$!} Then the chimera anti-diagonal
morphism could also be defined as: $z(c)=\Delta w\Rightarrow w\overset
{f(c)}{\longrightarrow}LimD$. There are always the over-and-back
factorizations of the ordinary intra-category maps at the top and bottom of
the adjunctive square (i.e., the upper and lower commutative triangles formed
by the anti-diagonal). The over-and-back factorization of the top map
$f(c)=z(c)h_{w}$ follows from the definition of the $z(c)$ just given:
\[
w\overset{h_{w}}{\Longrightarrow}\Delta w\overset{z(c)}{\Longrightarrow
}LimD=w\overset{h_{w}}{\Longrightarrow}\Delta w\Rightarrow w\overset
{f(c)}{\longrightarrow}LimD=w\overset{f(c)}{\longrightarrow}LimD.
\]

\noindent It was previously noted that some category theorists refer to the
natural transformation $\Delta w\overset{g(c)}{\longrightarrow}D$ as a
\textquotedblleft cone\textquotedblright\ while others refer to the
set-to-functor bundle of compatible maps---which is our chimera cone
$w\overset{c}{\Rightarrow}D$---as a cone. Often writers will effortlessly
switch back and forth between the two notions of a cone. The two inverse
chimeras $w\overset{h_{w}}{\Longrightarrow}\Delta w$ and $\Delta w\Rightarrow
w$ are used implicitly to switch back and forth, i.e., $\Delta w\Rightarrow
w\overset{c}{\Rightarrow}D=\Delta w\overset{g(c)}{\longrightarrow}D$ and
$w\overset{h_{w}}{\Longrightarrow}\Delta w\overset{g(c)}{\longrightarrow
}D=w\overset{c}{\Rightarrow}D$, which is the back and forth action of the
isomorphism for the representation on the left: $Hom(\Delta w,D)\cong
Het(w,D)$. That is used in the other over-and-back factorization
$g(c)=e_{D}z(c)$:

\begin{center}
$\Delta w\overset{z(c)}{\Longrightarrow}LimD\overset{e_{D}}{\Longrightarrow
}D=\Delta w\Rightarrow w\overset{f(c)}{\longrightarrow}LimD\overset{e_{D}%
}{\Rightarrow}D$

$=\Delta w\Rightarrow w\overset{c}{\Rightarrow}D=\Delta w\overset
{g(c)}{\longrightarrow}D$.
\end{center}

The remaining isomorphism $Het(w,D)\cong\mathcal{Z}(\Delta w,LimD)$ is the
isomorphism between cones $c:w\Rightarrow D$ and such cocones $z(c):\Delta
w\Rightarrow LimD$ (as defined above). The chimera zig-zag factorization is:%

\[
w\overset{c}{\Rightarrow}D=w\overset{h_{w}}{\Rightarrow}\Delta w\overset
{z(c)}{\Rightarrow}LimD\overset{e_{D}}{\Rightarrow}D.
\]

There are also the two universal anti-diagonal chimera morphisms. The cocone
$z(h_{w})=h_{w2}:\Delta w\Rightarrow Lim\Delta w$ has the diagonal map
$w\rightarrow Lim\Delta w$ as each component, and the cocone $z(e_{D}%
)=e_{D1}:\Delta LimD\Rightarrow LimD$ has the identity $1_{LimD}$ as each
component. By their universality properties, given any cocone of the form
$z(c):\Delta w\Rightarrow LimD$, there are unique morphisms $g=g(c):\Delta
w\rightarrow D$ and $f=f(c):w\rightarrow LimD$ (where $g=f^{\ast}$) that give
two factorizations of the anti-diagonal morphism:
\[
\Delta w\overset{h_{w2}}{\Longrightarrow}Lim\Delta w\overset{Limg}%
{\longrightarrow}LimD=\Delta w\overset{z(c)}{\Longrightarrow}LimD=\Delta
w\overset{\Delta f}{\longrightarrow}\Delta LimD\overset{e_{D1}}%
{\Longrightarrow}LimD.
\]

\noindent These two factorizations fit together to form the adjunctive-image
square with $z(c)$ as its main diagonal. Postcomposing the first equation with
$e_{D}$ yields the over-across-and-back factorization of $g(c):\Delta
w\rightarrow D$:

\begin{center}
$\Delta w\overset{h_{w2}}{\Longrightarrow}Lim\Delta w\overset{Limg}%
{\longrightarrow}LimD\overset{e_{D}}{\Longrightarrow}D=\Delta w\overset
{z(c)}{\Longrightarrow}LimD\overset{e_{D}}{\Longrightarrow}D=\Delta
w\overset{g(c)}{\longrightarrow}D$.
\end{center}

\noindent Precomposing the second equation with $h_{w}$ yields the
over-across-and-back factorization of $f(c):w\rightarrow LimD$:

\begin{center}
$w\overset{h_{w}}{\Rightarrow}\Delta w\overset{\Delta f}{\longrightarrow
}\Delta LimD\overset{e_{D1}}{\Longrightarrow}LimD=w\overset{h_{w}}%
{\Rightarrow}\Delta w\overset{z(c)}{\Longrightarrow}LimD=w\overset
{f(c)}{\longrightarrow}LimD$.
\end{center}

The two anti-diagonal universals also give the two over-and-back identities on
the functorial images (the `short forms' of the triangular identities):%

\begin{align*}
LimD\overset{h_{LimD}}{\Longrightarrow}\Delta LimD\overset{e_{D1}%
}{\Longrightarrow}LimD  &  =LimD\overset{1_{LimD}}{\longrightarrow}LimD\\
\Delta w\overset{h_{w2}}{\Longrightarrow}Lim\Delta w\overset{e_{\Delta w}%
}{\Longrightarrow}\Delta w  &  =\Delta w\overset{1_{\Delta w}}{\longrightarrow
}\Delta w.
\end{align*}

\section{Colimits in Sets}

Colimits in sets generalize the previous example of the coproduct
construction. We will consider colimits in $\mathbf{Set}$ but the argument
here (as with limits) would work for any other cocomplete (or complete)
category of algebras replacing the category of sets. The diagonal functor
$\Delta:\mathbf{Set}\rightarrow\mathbf{Set}^{\mathbf{D}}$ also has a
\textit{left} adjoint $Colim:\mathbf{Set}^{\mathbf{D}}\rightarrow\mathbf{Set}$.

For any diagram functor $D$ and set $z$, the adjunction for the colimit and
diagonal functors is:

\begin{center}
$Hom(ColimD,z)\cong Hom(D,\Delta z)$.
\end{center}

\noindent A adjunctive square for this adjunction has the form:

\begin{center}
$%
\begin{array}
[c]{ccc}%
(D,ColimD) & \overset{(f,Colimf)}{\longrightarrow} & (\Delta z,Colim\Delta
z)\\
\downarrow & \overset{^{(Colimf,\Delta g)}}{\nearrow} & \downarrow\\
(\Delta ColimD,ColimD) & \overset{(\Delta g,g)}{\longrightarrow} & (\Delta
z,z)
\end{array}
$
\end{center}

\noindent which commutes when $g=f^{\ast}$ or $f=g^{\ast}$.

For this adjunction, a chimera morphism or heteromorphism from a diagram
functor $D$ to a set $z$ is concretely given by a \textit{cocone}
$D\Rightarrow z$ which is defined as a set of maps $\{D_{i}\overset{g_{i}%
}{\longrightarrow}z\}$ indexed by the objects $i$ in the diagram category
$\mathbf{D}$ such that for any morphism $\alpha:i\rightarrow j$ in
$\mathbf{D}$, $D_{i}\overset{D_{\alpha}}{\longrightarrow}D_{j}\overset{g_{j}%
}{\longrightarrow}z=D_{i}\overset{g_{i}}{\longrightarrow}z$. The adjunction is
then given by the birepresentations of the het-bifunctor where
$Het(D,z)=\{D\Rightarrow z\}$ is the set of cocones from the diagram functor
$D$ to the set $z$.

Since the role of the diagram and set are reversed from the case of limits, a
new notion of \textquotedblleft coelements\textquotedblright\ as
\textquotedblleft determiners\textquotedblright\ from $D$ is necessary.
Conceptually, if we go back to the idea of a function as a way for
(co)elements in the domain to determine certain elements in the codomain, then
a \textit{coelement\ of }$D$ is the minimal set in the domain $D_{i}$s that
are necessary for functions on those domains to compatibly determine a one
point set as the codomain. This could be thought of as a minimal
partially-defined 1-cocone, or better, the germ of a cocone to determine a one
point set before such a set is selected. If the diagram category was discrete
so there were no maps between the $D_{i}$s, then a coelement would simply be a
member $x_{i}$ of one of the $D_{i}$s since a function defined on that point
to the one point set is sufficient to \textquotedblleft
determine\textquotedblright\ that point in the codomain. The set of coelements
of $D$ would then be the coproduct (disjoint union) of the $D_{i}$s. But if
there were maps between the $D_{i}$s such as $D_{\alpha}:D_{i}\rightarrow
D_{j}$ then $x_{j}=D_{\alpha}(x_{i})$ would also need to be mapped to the same
one point. Hence we define a compatibility equivalence relation on the
disjoint union of the $D_{i}$s where $x_{i}\sim x_{j}$ if $D_{\alpha}%
(x_{i})=x_{j}$ for any morphism $D_{\alpha}$ between the $D_{i}$s. Thus a
\textit{coelement} or atomic determiner (\textquotedblleft germ of a
cocone\textquotedblright) from $D$ would consist of an equivalence class or
block in the partition of the disjoint union determined by the compatibility
equivalence relation.

Each coelement of $D$ represents a determiner (without a specified determinee
point), a germ of functions on the $D_{i}$s as domains to compatibly determine
a single point in the codomain. As always, the left adjoint repackages an
object $D$ in the top category as the object of determiners $ColimD$ in the
bottom category so that a determination would be represented by a morphism in
the bottom category. Thus a chimera $D\Rightarrow z$ (i.e., a cocone) would be
represented by a set morphism $g:ColimD\rightarrow z$.

As before (but dually), the basic conceptual move is to take such a conceptual
atom, a determiner coelement of $D$, to be its own determinee, and that
self-determination is represented by the identity map $1_{ColimD}%
:ColimD\rightarrow ColimD$, the set of determiners reconceptualized as the
elements of the set of determinees. But for that self-determination from $D$
to be represented by a morphism in the top category with domain $D$, the right
adjoint must, as always, repackage the determinees of an object in the bottom
category as an object in the top category. The right adjoint $\Delta$
repackages the determinees or coelements of the set $ColimD$ as the
determinees of the functor $\Delta ColimD$ and the determination from $D$ is
expressed by the canonical morphism $\eta_{D}:D\rightarrow\Delta ColimD$.

This self-determination is much clearer in the one-functor treatment of the
chimera version of the unit. The chimera unit $h_{D}:D\Rightarrow ColimD$ is
simply the cocone of injection maps. The coelements of $D$ are the determiners
so if they are collected together as the determinees in a set $ColimD$, then
the self-determination would be expressed by the heteromorphic unit
$h_{D}:D\Rightarrow ColimD$, namely the injection maps which map each
coelement of $D$ to itself as an element of $ColimD$.

If the diagram category $\mathbf{D}$ is multiply connected, then the
coelements (determiners or germs) of any diagram $D$ would correspond to the
coproduct of the coelements of $D$ restricted to each component. In the two
component case, the coelements of $D$ would correspond to the members of the
coproduct or disjoint sum $ColimD\mathbf{\upharpoonright}%
1+ColimD\mathbf{\upharpoonright}2$.

Starting with a set $z$, its determinees are simply its elements as a set.
\ The right adjoint $\Delta z$ represents those elements as the determinees of
a diagram functor so that a chimera morphism $D\Rightarrow z$ would be
expressed by a morphism in the top category (natural transformation)
$f:D\rightarrow\Delta z$. \ But as a diagram functor, $\Delta z$ can also be a
determiner in a determination to $z$ by the identity morphism $1_{\Delta
z}:\Delta z\rightarrow\Delta z$. For that self-determination to be expressed
by a morphism in the bottom category (i.e., a set map), the determiners of
$\Delta z$ must be repackaged by the left adjoint as determiners $Colim\Delta
z$ in the bottom category and that self-determination is realized by the
canonical map $\varepsilon_{z}:Colim\Delta z\rightarrow z$. The colimit
$Colim\Delta z$ is the coproduct of a copy of $z$, one for each connected
component in the diagram category. \ Thus for the case of two components, the
counit $\varepsilon_{z} $ is the codiagonal or folding map $Colim\Delta z\cong
z+z\rightarrow z$. The chimera version $e_{z}:\Delta z\Rightarrow z$ is simply
the cocone of identity maps $1_{z}$. The determinees of $z$ are repackaged as
the determiner coelement of $\Delta z$ and then the heteromorphic counit
$e_{z}:\Delta z\Rightarrow z$ gives the self-determination wherein the
coelements of $\Delta z$ determine themselves as elements of $z$.

Since the limit and colimit are respectively the right and left adjoints to
the same constant functor, the two adjunctions give bidirectional
determination from sets $w$ to diagram functors $D$ (e.g., cones) and from
diagram functors $D$ to sets $z$ (e.g., cocones).

Here again, the concrete heteromorphisms can be used to give a chimera version
of the adjunctive squares diagram. Given a cocone $c:D\Rightarrow z$, there is
a unique natural transformation $f(c):D\rightarrow\Delta z$ such that
$D\overset{c}{\Rightarrow}z=D\overset{f(c)}{\rightarrow}\Delta z\overset
{e_{z}}{\Rightarrow}z$ which gives the representation on the right:
$Het(D,z)\cong Hom(D,\Delta z)$. And there is a unique set map
$g(c):ColimD\rightarrow z$ such that $D\overset{c}{\Rightarrow}z=D\overset
{h_{D}}{\Rightarrow}ColimD\overset{g(c)}{\rightarrow}z$ which gives the
representation on the left: $Hom(ColimD,z)\cong Het(D,z)$. Splicing the
commutative triangles together gives the adjunctive square with $D\overset
{c}{\Rightarrow}z$ as the main diagonal.

\begin{center}
$%
\begin{array}
[c]{ccccc}
& D & \overset{f(c)}{\longrightarrow} & \Delta z & \\
h_{D} & \Downarrow &  & \Downarrow & e_{z}\\
& ColimD & \overset{g(c)}{\longrightarrow} & z &
\end{array}
\medskip$

Chimera Adjunctive Square for Colimit Adjunction
\end{center}

\noindent The anti-diagonal map $z(c):ColimD\Rightarrow\Delta z$ is the
set-to-functor cone each of whose components is $g(c):ColimD\rightarrow z.$ It
could also be constructed using the chimera cone $z\Rightarrow\Delta z$ (each
component is $1_{z}$ so it is the chimera universal denoted $h_{w}%
:w\Rightarrow\Delta w$ in the adjunction for limits) that is inverse to the
cocone $e_{z}:\Delta z\Rightarrow z$. Then $z(c)=ColimD\overset{g(c)}%
{\longrightarrow}z\Rightarrow\Delta z$ which also establishes the remaining
isomorphism: $Hom(ColimD,z)\cong\mathcal{Z}(ColimD,\Delta z)$. The
over-and-back factorization of $g(c)$ is immediate:

\begin{center}
$ColimD\overset{z(c)}{\Longrightarrow}\Delta z\overset{e_{z}}{\Longrightarrow
}z=ColimD\overset{g(c)}{\longrightarrow}z\Rightarrow\Delta z\overset{e_{z}%
}{\Rightarrow}z=ColimD\overset{g(c)}{\longrightarrow}z$.
\end{center}

\noindent The other over-and-back factorization is:

\begin{center}
$D\overset{h_{D}}{\Rightarrow}ColimD\overset{z(c)}{\Rightarrow}\Delta
z=D\overset{h_{D}}{\Rightarrow}ColimD\overset{g(c)}{\rightarrow}%
z\Rightarrow\Delta z$

$=D\overset{c}{\Rightarrow}z\Rightarrow\Delta z=D\overset{f(c)}{\rightarrow
}\Delta z$
\end{center}

\noindent where the last equality is the way of going from cocones as chimera
to \textquotedblleft cocones\textquotedblright\ as natural transformations.

The isomorphisms of the adjunction have been established:

\begin{center}
$Hom(ColimD,z)\cong Het(D,z)\cong\mathcal{Z}(ColimD,\Delta z)\cong
Hom(D,\Delta z)$.
\end{center}

\noindent The zig-zag factorization is:%

\[
D\overset{c}{\Rightarrow}z=D\overset{h_{D}}{\Rightarrow}ColimD\overset
{z(c)}{\Rightarrow}\Delta z\overset{e_{z}}{\Rightarrow}z.
\]

As in the case of the limits, there are the two anti-diagonal chimera
universals with the analogous properties.

\section{Adjoints to Forgetful Functors}

Perhaps the most accessible adjunctions are the free-forgetful adjunctions
between $\mathbf{X}=\mathbf{Set}$ and a category of algebras such as the
category of groups $\mathbf{A}=\mathbf{Grps}$. \text{The right adjoint
}$G:\mathbf{A}\rightarrow\mathbf{X}$ forgets the group structure to give the
underlying set $Ga$ of a group $a$. \ The left adjoint $F:\mathbf{X}%
\rightarrow\mathbf{A}$ gives the free group $Fx$ generated by a set $x$. The
hom-set isomorphism and the adjunctive squares have the usual forms.

For this adjunction, the heteromorphisms are the set-to-group functions
$x\overset{c}{\Rightarrow}a$ and the het-bifunctor is given by such functions:
$Het(x,a)=\{x\Rightarrow a\}$. A chimera $c:x\Rightarrow a$ determines a set
map $f=f(c):x\rightarrow Ga$ and a group homomorphism $g(c)=f^{\ast
}:Fx\rightarrow a$ so that $\widehat{x}=(x,Fx)\overset{(f,f^{\ast}%
)}{\longrightarrow}(Ga,a)=\widehat{a}$ is the abstract version of the concrete
$x\overset{c}{\Rightarrow}a$. These associations also give us the two representations:%

\[
Hom(Fx,a)\cong Het(x,a)\cong Hom(x,Ga).
\]

\noindent The universal element for the functor $Het(x,-)$ is the chimera
$h_{x}:x\Rightarrow Fx$ (insertion of the generators into the free group) and
the universal element for the functor $Het(-,a)$ is the chimera $e_{a}%
:Ga\Rightarrow a$ (the retracting of the elements of the underlying set back
to the group).

The right adjoint always gives a representation of all the possible
determinees of the target $a$ as an $\mathbf{X}$-object $Ga$ with maximal
structure so that a determination $x\Rightarrow a$ would be represented by an
$\mathbf{X}$-morphism $x\rightarrow Ga$. The underlying set functor does
exactly that and a set map $x\rightarrow Ga$ gives such a determination. \ The
left adjoint always gives a representation of the determiners of a source $x$
as an $\mathbf{A}$-object $Fx$ with minimal structure so that a determination
$x\Rightarrow a$ would be represented by an $\mathbf{A}$-morphism
$Fx\rightarrow a$. Since $\mathbf{A}$ is the category of groups, the
representation of the elements of a set $x$ as elements of a group independent
of any target is accomplished by the free group $Fx$. No extra
\textquotedblleft junk\textquotedblright\ is added other than the group
elements generated by the generators $x$ and there is no
noise\ (\textquotedblleft noise\textquotedblright\ in the sense of identifying
distinct \textquotedblleft signals\textquotedblright) and no extra relations
are imposed (other than those necessary to make it a group). \ Then the image
of the generators under a group homomorphism $Fx\rightarrow a$ is a
determination from the set elements of $x$ to the group elements of $a$, and
every such determination would generate such a group homomorphism.

The set-to-group self-determination $Ga$-to-$a$ is correlated with the set map
$1_{Ga}:Ga\rightarrow Ga$. For the determiners from a set to be represented as
an $\mathbf{A}$-object, the left adjoint $F$ must be applied so, in this case,
the free group functor yields the free group $FGa$ generated by the underlying
set of the group $a$. \ Thus $FGa$ is the determinees $Ga$ of $a$ represented
as determiners and the group homomorphism induced by the set map $1_{Ga}$
gives that self-determination as the counit $\varepsilon_{a}:FGa\rightarrow a$
(which is the adjoint transpose of $1_{Ga}$). The simpler chimera counit
$e_{a}:Ga\Rightarrow a$ just retracts the elements of the underlying set back
to the group. The determinee elements of the group $a$ are represented as the
determinee elements of the set $Ga$ but then they turn around and determine
themselves by the chimera counit $e_{a}:Ga\Rightarrow a$ (which is the adjoint
correlate of $1_{Ga}$).

The set-to-group self-determination $x$-to-$Fx$ is correlated with the group
homomorphism $1_{Fx}:Fx\rightarrow Fx$. For the determinees of a group to be
represented as an $\mathbf{X}$-object, the right adjoint $G$ must be applied
so, in this case, the forgetful functor yields the underlying set $GFx$ of the
free group generated by the set $x$. Thus $GFx$ is the determiners $Fx$ from
$x$ represented as determinees and the set map induced by the homomorphism
$1_{Fx}$ gives that self-determination as the unit $\eta_{x}:x\rightarrow GFx$
(the adjoint transpose of $1_{Fx}$) The simpler chimera version $h_{x}%
:x\Rightarrow Fx$ is the injection of the generators into the free group. The
determiner elements of the set $x$ are represented as determiner elements of
the free group $Fx$ but then they are determined by themselves via the chimera
unit $h_{x}:x\Rightarrow Fx$ (the adjoint correlate of $1_{Fx}$).

In an adjunctive square, any set-to-group determination $(f,f^{\ast})$
expressed by a set map $f:x\rightarrow Ga$ and its associated group
homomorphism $f^{\ast}:Fx\rightarrow a$ would be factored through both the
sending universal $(\eta_{x},1_{Fx})$ and the receiving universal
$(1_{Ga},\varepsilon_{a})$ by the indirect anti-diagonal map $(Gf^{\ast},Ff) $.

The chimera version of the adjunctive squares diagram always has the generic
form where $x\overset{c}{\Rightarrow}a$ is the main diagonal.

\begin{center}
$%
\begin{array}
[c]{ccccc}
& x & \overset{f(c)}{\longrightarrow} & Ga & \\
h_{x} & \Downarrow &  & \Downarrow & e_{a}\\
& Fx & \overset{g(c)}{\longrightarrow} & a &
\end{array}
$
\end{center}

\noindent In this case, the anti-diagonal map $z(c):Fx\Rightarrow Ga$ is
essentially the group homomorphism $g(c):Fx\rightarrow a$ but where the
codomain is taken as the underlying set. Intuitively, there is an inverse
$s_{a}:a\Rightarrow Ga$ to $e_{a}:Ga\Rightarrow a$ such that $e_{a}s_{a}=1_{a}
$ and $s_{a}e_{a}=1_{Ga}$ so $Fx\overset{z(c)}{\Rightarrow}Ga=Fx\overset
{g(c)}{\rightarrow}a\overset{s_{a}}{\Rightarrow}Ga$. The proof of the
over-and-back factorization for $f:x\rightarrow Ga$ goes around the square counter-clockwise:

\begin{center}
$x\overset{h_{x}}{\Longrightarrow}Fx\overset{z(c)}{\Rightarrow}Ga=x\overset
{h_{x}}{\Longrightarrow}Fx\overset{g(c)}{\rightarrow}a\overset{s_{a}%
}{\Rightarrow}Ga=x\overset{c}{\Rightarrow}a\overset{s_{a}}{\Rightarrow
}Ga=x\overset{f(c)}{\longrightarrow}Ga\overset{e_{a}}{\Rightarrow}%
a\overset{s_{a}}{\Rightarrow}Ga=x\overset{f(c)}{\longrightarrow}Ga$
\end{center}

\noindent while the other over-and-back factorization for $g:Fx\rightarrow a$
is trivial:

\begin{center}
$Fx\overset{z(c)}{\Rightarrow}Ga\overset{e_{a}}{\Longrightarrow}%
a=Fx\overset{g(c)}{\rightarrow}a\overset{s_{a}}{\Rightarrow}Ga\overset{e_{a}%
}{\Longrightarrow}a=Fx\overset{g(c)}{\rightarrow}a$.
\end{center}

The chimera isomorphism between the set-to-group main diagonals and the
group-to-set anti-diagonals is: $Het(x,a)=\{x\overset{c}{\Rightarrow}%
a\}\cong\{Fx\overset{z(c)}{\Rightarrow}Ga\}=\mathcal{Z}(Fx,Ga)$. The zig-zag
factorization is:%

\[
x\overset{c}{\Rightarrow}a=x\overset{h_{x}}{\Rightarrow}Fx\overset
{z(c)}{\Rightarrow}Ga\overset{e_{a}}{\Rightarrow}a.
\]

The two anti-diagonal universals are $h_{x2}=z(h_{x}):Fx\Rightarrow GFx$ which
injects a free group into its underlying set, and $e_{a1}=z(e_{a}%
):FGa\Rightarrow Ga$ which canonically maps the free group on the underlying
set of a group $a$ to that underlying set $Ga$ using the group operations of
$a$. One over-and-back identity is: $Fx\overset{h_{x2}}{\Longrightarrow
}GFx\overset{e_{Fx}}{\Longrightarrow}Fx=Fx\overset{1_{Fx}}{\longrightarrow}%
Fx$. The map of a free group onto its underlying set and back to the group is
trivially the identity group homomorphism. And the other over-and-back
identity is: $Ga\overset{h_{Ga}}{\Longrightarrow}FGa\overset{e_{a1}%
}{\Longrightarrow}Ga=Ga\overset{1_{Ga}}{\longrightarrow}Ga$. The injection of
the underlying set $Ga$ of a group into the free group on that set followed by
the projection of that free group onto $Ga$ using the group operations of $a$
is the identity on $Ga$. There are also the two over-across-and-back
factorizations that might be checked:%

\begin{align*}
x\overset{h_{x}}{\Longrightarrow}Fx\overset{Ff}{\longrightarrow}%
FGa\overset{e_{a1}}{\Rightarrow}Ga  &  =x\overset{f}{\longrightarrow}Ga\\
Fx\overset{h_{x2}}{\Longrightarrow}GFx\overset{Gg}{\longrightarrow}%
Ga\overset{e_{a}}{\Longrightarrow}a  &  =Fx\overset{g}{\longrightarrow}a.
\end{align*}

The four naturally isomorphic bifunctors are:%

\[
Hom(Fx,a)\cong Het(x,a)\cong\mathcal{Z}(Fx,Ga)\cong Hom(x,Ga)
\]

\noindent and the four associated maps in the same order are:%

\[
Fx\overset{g(c)}{\longrightarrow}a\leftrightsquigarrow x\overset
{c}{\Longrightarrow}a\leftrightsquigarrow Fx\overset{z(c)}{\Longrightarrow
}Ga\leftrightsquigarrow x\overset{f(c)}{\longrightarrow}Ga.
\]

\noindent This includes all the four possibilities of group-to-group,
set-to-group, group-to-set, and set-to-set morphisms. The mixed or mongrel
morphisms are the heteromorphisms while the unmixed or pure morphisms are the
homomorphisms in the standard theory. Note the asymmetry in that set-to-group
determination is for arbitrary sets $x$ and arbitrary groups $a$ whereas the
chimera morphisms in the opposite direction are only defined on the images of
the two functors, i.e., are only from free groups to underlying sets of
groups. This adjunction does not contemplate arbitrary determination from
groups to sets. If the underlying set functor had a right adjoint, then there
would be arbitrary two-way determination between sets and groups expressed by
adjunctions. \ But it does not have a right adjoint for reasons that will now
be explained.

Adjoints to underlying set functors $U$ are particularly accessible since the
underlying determinations can be expressed by easily understood chimera
functions between the objects in the two categories. \ Given such a function
$x\Rightarrow a$\ (where $a$ is from the category of objects with structure
and $x$ is a set), the existence of a left adjoint to $U$ (i.e., a left
representation of $Het(x,a)=\{x\Rightarrow\alpha\}$) will depend on whether or
not there is an $\mathbf{A}$-object $Fx$ with the least or minimal structure
so that every determination $x\rightarrow Ua$ by a set morphism will have an
adjoint transpose by an $\mathbf{A}$\textbf{-}morphism $Fx\rightarrow a$. But
there might also be a chimera $a\Rightarrow x$\ where the existence of a
\textit{right} adjoint to $U$ will depend on whether or not there is an
$\mathbf{A}$-object $Ix$ with the greatest or maximum structure so that any
set function $Ua\rightarrow x$ can also be expressed by an $\mathbf{A}%
$-morphism $a\rightarrow Ix$. The homomorphisms $a\rightarrow Ix$ would have
to preserve the structure which was forgotten in $Ua\rightarrow x$ so $Ix$
would have to carry all the possible structures that might carried by the
structure-preserving morphisms $a\rightarrow Ix$. There is no `maximal group'
$Ix$ so that $a\rightarrow Ix$ could be the adjoint transpose to
$Ua\rightarrow x$; hence the functor giving the underlying set of a group has
no right adjoint. \ 

For another example, consider the underlying set functor $U:\mathbf{Pos}%
\rightarrow\mathbf{Set}$ from the category of partially ordered sets (an
ordering that is reflexive, transitive, and anti-symmetric) with
order-preserving maps to the category of sets. It has a left adjoint since
each set has a least partial order on it, namely the discrete ordering.
\ Hence any chimera function $x\Rightarrow a$\ from a set $x$ to a partially
ordered set or poset $a$ could be expressed as a set function $x\rightarrow
Ua$ or as an order-preserving function $Dx\rightarrow a$ where $Dx$ gives the
discrete ordering on $x$. In the other direction, one can have a chimera
function $a\Rightarrow x$\ and an ordinary set function $Ua\rightarrow x$ but
the underlying set functor $U$ does not have a right adjoint since there is no
maximal partial order $Ix$ on $x$ so that any determination $Ua\rightarrow x$
could be expressed as an order-preserving function $a\rightarrow Ix$. To
receive all the possible orderings, the ordering relation would have to go
both ways between any two points which would then be identified by the
anti-symmetry condition so that $Ix$ would collapse to a single point. Thus
poset-to-set determinations expressed by $a\Rightarrow x $\ cannot be
represented as determination through universals, i.e., by an adjunction.

Relaxing the anti-symmetry condition, let $U:\mathbf{Ord}\rightarrow
\mathbf{Set}$ be the underlying set functor from the category of preordered
sets (reflexive and transitive orderings) to the category of sets. \ The
discrete ordering again gives a left adjoint. \ But now there is also a
maximal ordering on a set $x$, namely the `indiscrete' ordering $Ix$ on $x$
(the `indiscriminate' or `chaotic' preorder on $x$) which has the ordering
relation both ways between any two points. \ Then a preorder-to-set chimera
morphism $a\Rightarrow x$\ (just a set function ignoring the ordering) can be
represented either as a set function $Ua\rightarrow x$ or as an
order-preserving function $a\rightarrow Ix$ so that $U$ also has a right
adjoint $I$. Thus the determination both ways between preordered sets and sets
can be given by adjunctions, i.e., represented as determination through universals.

\section{The Product-Exponential Adjunction in Sets}

The product-exponential adjunction in $\mathbf{Set}$ is an interesting example
of an adjoint pair of endo-functors $F:\mathbf{Set}\rightleftarrows
\mathbf{Set}:G$. For any fixed (non-empty) \textquotedblleft
index\textquotedblright\ set $A$, the product functor $F(-)=-\times
A:\mathbf{Set}\rightarrow\mathbf{Set}$ has a right adjoint $G(-)=(-)^{A}%
:\mathbf{Set}\rightarrow\mathbf{Set}$ which makes $\mathbf{Set}$ a cartesian
closed category. For any sets $X$ and $Y$, the adjunction has the form:

\begin{center}
$Hom(X\times A,Y)\cong Hom(X,Y^{A})$.
\end{center}

Since both functors are endo-functors on $\mathbf{Set}$, we cannot expect to
find any chimera morphisms outside this category. The job of finding some
`chimera' morphisms of $\mathbf{Set}$ such that the adjunction arises out of
birepresenting them is trivial; take either of the hom-sets $Hom(X\times A,Y)
$ or $Hom(X,Y^{A})$. But that does not show that this theory of adjoints works
for the product-exponential adjunction since we don't have the chimera
universals. For instance, the chimera sending universal should be a canonical
morphism $h_{X}:X\Rightarrow FX$ but if $FX=X\times A$, there is no canonical
map $X\rightarrow X\times A$ (except in the special case where $A$ is a
singleton). Similarly, the receiving universal should be a canonical map
$e_{Y}:GY\Rightarrow Y$ but if $GY=Y^{A}$ then there is no canonical map
$Y^{A}\rightarrow Y$ (unless $A$ is a singleton). Hence there appears to be a
problem in applying the theory of adjoints to the product-exponential
adjunction (or to any adjunction where one category is a subcategory of the other).

The way out of this apparent problem is shown by the following general result
\cite[p. 83]{Freyd:ac}. Consider any adjunction $F:\mathbf{X}\rightleftarrows
\mathbf{A}:G$. Let $X^{\prime}$\textbf{\ }be the subcategory of $\mathbf{X}$
generated by the image of $G$ where if $G$ is one-one on objects, then its
image is that subcategory. Dually, let $\mathbf{A}^{\prime}$ be the
subcategory of $\mathbf{A}$ generated by $F$ (where the image is the
subcategory if $F$ is one-one on objects). If $x^{\prime}=Ga$, then

\begin{center}
$Hom_{\mathbf{X}^{\prime}}(GFx,x^{\prime})\cong Hom_{\mathbf{X}^{\prime}%
}(GFx,Ga)\cong Hom_{\mathbf{A}}(Fx,a)\cong Hom_{\mathbf{X}}(x,Ga)\cong
Hom_{\mathbf{X}}(x,x^{\prime})$.
\end{center}

\noindent Thus $\mathbf{X}^{\prime}$ is a reflective subcategory of
$\mathbf{X}$ where the reflector (left adjoint to the inclusion) can be taken
as the functor $GF:\mathbf{X}\rightarrow\mathbf{X}^{\prime}$ (where $G$ was
construed as taking values in $\mathbf{X}^{\prime}$).

Dually, let $a^{\prime}=Fx$ so that

\begin{center}
$Hom_{\mathbf{A}}(a^{\prime},a)\cong Hom_{\mathbf{A}}(Fx,a)\cong
Hom_{\mathbf{X}}(x,Ga)\cong Hom_{\mathbf{A}}(Fx,FGa)\cong Hom_{\mathbf{A}%
^{\prime}}(a^{\prime},FGa)$.
\end{center}

\noindent Thus $\mathbf{A}^{\prime}$ is a coreflective subcategory of
$\mathbf{A}$ where the coreflector (right adjoint to the inclusion) can be
taken as the functor $FG:\mathbf{A}\rightarrow\mathbf{A}^{\prime}$ (where $F$
is construed as taking values in $\mathbf{A}^{\prime}$).

It was previously noted that in the case of a reflection, i.e., a left adjoint
to the inclusion functor, heteromorphisms can be found as the morphisms with
their tail in the ambient category and their heads in the subcategory. For a
coreflection (right adjoint to the inclusion functor), the heteromorphisms
would be turned around, i.e., would have their tail in the subcategory and
their head in the ambient category. If neither $\mathbf{X}$ nor $\mathbf{A}$
were a subcategory of the other, then the above construction would not find
the true chimera morphisms from objects in $\mathbf{X}$ to objects in
$\mathbf{A}$. It would only find what might be viewed as \textquotedblleft
pseudo-chimera\textquotedblright\ morphisms since in the reflective case, a
morphism $x\rightarrow x^{\prime}$ is just a morphism $x\rightarrow Ga$ and
$Ga$ is not an $\mathbf{A}$-object at all. Or in the coreflective case, a
morphism $a^{\prime}\rightarrow a$ is only a pseudo-chimera morphism of the
form $Fx\rightarrow a$ since $Fx$ is not an $\mathbf{X}$-object at all. But if
$\mathbf{X}$ or $\mathbf{A}$ is a subcategory of the other (including the case
$\mathbf{X}=\mathbf{A}$), then there is nothing \textquotedblleft
pseudo\textquotedblright\ about the chimeras identified in the above
reflective and coreflective cases. \textit{Given} that one category is a
subcategory of the other, that is as \textquotedblleft
hetero\textquotedblright\ as the chimeras can be--and that is the case at
hand. In the extreme case of the ur-adjunction (the self-adjunction of the
identity functor on any category), all differences between the homo- and
heteromorphisms are wiped out.

In this case, $\mathbf{X}=\mathbf{Set}=\mathbf{A}$. There are dual ways of
reinterpreting the adjunction---either as a reflection or coreflection. As a
reflection, let $\mathbf{APower}$ be the subcategory of $G(-)=(-)^{A}$ images
($G$ is one-one since $A$ is non-empty) so that $\mathbf{APower}%
\hookrightarrow\mathbf{Set}$, and that inclusion functor has a left adjoint
$F^{\prime}(-)=(-\times A)^{A}:\mathbf{Set}\rightarrow\mathbf{APower}$. Then
the heteromorphisms are those with their tail in $\mathbf{Set}$ and head in
$\mathbf{APower}$, i.e., the morphisms of the form $X\rightarrow Y^{A}$. But
now we have the chimera universals. The chimera unit $h_{X}:X\Rightarrow
F^{\prime}X=(X\times A)^{A}$ is the canonical map that takes an $x$ in $X$ to
the function $(x,-):A\rightarrow X\times A$ which takes $a$ in $A$ to
$(x,a)\in X\times A$ which is also the unit $\eta_{X}:X\rightarrow(X\times
A)^{A}$ in the original product-exponential adjunction. Since the right
adjoint in the reflective case is the inclusion, the chimera counit $e_{Y^{A}%
}:Y^{A}\rightarrow Y^{A}$ is the identity. The chimera adjunctive square then
is the following commutative diagram.

\begin{center}
$%
\begin{array}
[c]{ccccc}
& X & \overset{f}{\rightarrow} & Y^{A} & \\
h_{X} & \Downarrow &  & \Downarrow & e_{Y^{A}}\\
& (X\times A)^{A} & \overset{(f^{\ast})^{A}}{\longrightarrow} & Y^{A} &
\end{array}
$
\end{center}

As a coreflection, let $\mathbf{AProd}$ be the subcategory of $F(-)=-\times A
$ images ($F$ is one-one since $A$ is non-empty) so that $\mathbf{AProd}%
\hookrightarrow\mathbf{Set}$, and that inclusion functor has a right adjoint
$G^{\prime}(-)=(-)^{A}\times A:\mathbf{Set}\rightarrow\mathbf{AProd}$. Then
the heteromorphisms are those with their tail in $\mathbf{AProd}$ and their
head in $\mathbf{Set}$, i.e., the morphisms of the form $X\times A\rightarrow
Y$. And now we again have the chimera universals. The chimera counit
$e_{Y}:Y^{A}\times A=G^{\prime}Y\Rightarrow Y$ is the evaluation map which is
also the counit $\varepsilon_{Y}:Y^{A}\times A\rightarrow Y$ in the original
product-exponential adjunction. Since the left adjoint in the coreflective
case is the inclusion, the chimera unit $h_{X\times A}:X\times A\Rightarrow
X\times A$ is the identity. The chimera adjunctive square is then the
following commutative diagram.

\begin{center}
$%
\begin{array}
[c]{ccccc}
& X\times A & \overset{g^{\ast}\times A}{\longrightarrow} & Y^{A}\times A & \\
h_{X\times A} & \Downarrow &  & \Downarrow & e_{Y}\\
& X\times A & \overset{g}{\rightarrow} & Y &
\end{array}
$
\end{center}

Hence there are in fact \textit{two} ways of choosing the heteromorphisms and
each way determines chimera universals which have all the usual factorizations
and identities holding. A similar treatment would work for any other case of
an adjunction $F:\mathbf{X}\rightleftarrows\mathbf{A}:G$ where one of the
categories was a subcategory of the other.

\section{Galois Connections}

A partially ordered set or poset is construed as a category where each hom-set
either has one map (the relation $\leq$ holds) or no maps (the relation $\leq$
does not hold). A functor between posets is an order-preserving map. An
adjunction between posets is usually known as a \textit{Galois connection}
$\cite[p. 93]{mac:cwm}$.

A standard Galois connection is the direct image and inverse image adjunction
induced by any function $f:\mathcal{X}\rightarrow\mathcal{A}$ between sets
$\mathcal{X}$ and $\mathcal{A}$. Let $\mathbf{X}$ be the power set of a set
$\mathcal{X}$ and let $\mathbf{A}$ be the power set of $\mathcal{A}$ both with
the inclusion order. Then $F=f():\mathbf{X}\rightarrow\mathbf{A}$, the direct
image map, and $G=f^{-1}:\mathbf{A}\rightarrow\mathbf{X}$, the inverse image
map, are both order-preserving functions. The adjunction, $Fx\leq a$ iff
$x\leq Ga$, will be written as: $f(x)\subseteq a$ iff $x\subseteq f^{-1}(a)$
for any subsets $x$ and $a$.

This adjunction is about the determination a subset $x$ of $\mathcal{X}$ to a
subset $a$ of $\mathcal{A}$ by using the function $f:\mathcal{X}%
\rightarrow\mathcal{A}$ where the relation $x\Rightarrow a$ holds if for all
$\zeta\in x$, $f(\zeta)\in a$ (which is just a point-wise way of saying that
the direct image $f(x)$ is a subset of $a$). A value $Het(x,a)=\{x\Rightarrow
a\}$ is a singleton if $x\Rightarrow a$ holds and is empty otherwise.

As always, the right adjoint applied to an $\mathbf{A}$-object $a$ gives the
$\mathbf{X}$-object $Ga$ that represents all the possible determinees of $a$.
The determinees of $a$ are its subsets of the form $f(x)$ for some $x$ in
$\mathbf{X}$ so the representation of all the possible determinees \textit{as
an }$\mathbf{X}$\textit{-object} would obtained as the union $\cup\{x\mid
f(x)\subseteq a\}=f^{-1}(a)$ or in general for Galois connections,
$\sup\{x\mid Fx\leq a\}=Ga$.

As always, the left adjoint applied to an $\mathbf{X}$-object $x$ gives the
$\mathbf{A}$-object $Fx$ that represents all the possible determiners from
$x$. A set $x$ determines a subset of an $\mathbf{A}$-object $a$ whenever
$x\subseteq f^{-1}(a)$ so the representation of all the possible
determinations from $x$ as an $\mathbf{A}$-object could be obtained as
$\cap\{a\mid x\subseteq f^{-1}(a)\}=f(x)$ or in general for Galois
connections, $\inf\{a\mid x\leq Ga\}=Fx$.

As indicated by the inclusion $f^{-1}(a)\subseteq f^{-1}(a)$, the $\mathbf{X}
$-object that represents the determinees of $a$ can also be a source of a
subset of $a$. The representation of $f^{-1}(a)$ as an $\mathbf{A}$-object
determiner is obtained by applying the left adjoint so we have the canonical
inclusion: $f\left(  f^{-1}(a)\right)  \subseteq a$ which also gives the true
chimera relation $e_{a}:f^{-1}(a)\Rightarrow a$ between subsets of different sets.

As indicated by the inclusion $f\left(  x\right)  \subseteq f(x)$, the
$\mathbf{A}$-object that represents the determiners from $x$ can also be a
target of $x$'s determination. \ The application of the right adjoint gives
the $\mathbf{X}$-object $f^{-1}(f(x))$ so that the determination from $x$
takes the form of the canonical inclusion: $x\subseteq f^{-1}(f(x))$ which
also gives the true chimera relation $h_{x}:x\Rightarrow f(x)$ between subsets
of different sets.

The chimera version of the adjunctive square then has the following form of an
\textquotedblleft if and only if\textquotedblright\ (iff) statement.

\begin{center}
$%
\begin{array}
[c]{ccccc}
& x & \subseteq & f^{-1}(a) & \\
h_{x} & \Downarrow &  & \Downarrow & e_{a}\\
& f(x) & \subseteq & a &
\end{array}
$
\end{center}

\noindent Since the vertical relations always hold, the statement is:
$f(x)\subseteq a$ iff $x\Rightarrow a$ iff $x\subseteq f^{-1}(a)$ (the
representation isomorphisms in this case). The anti-diagonal relation could be
defined by: $f(x)\Rightarrow f^{-1}(a)$ if $f^{-1}(f(x))\subseteq f^{-1}(a)$
or equivalently as: $f(x)\subseteq f(f^{-1}(a))$. One over-and-back
factorization is $x\subseteq f^{-1}(a)$ iff $x\subseteq f^{-1}(f(x)\subseteq
f^{-1}(a)$, and the other one is: $f(x)\subseteq a$ iff $f(x)\subseteq
f(f^{-1}(a))\subseteq a$. Then the zig-zag factorization is the statement:
$x\Rightarrow a$ iff $x\Rightarrow f(x)\Rightarrow f^{-1}(a)\Rightarrow a$, or
in conventional terms: $f(x)\subseteq a$ iff $f(x)\subseteq f(x)\subseteq
f(f^{-1}(a))\subseteq a$. The chimera `isomorphism' is: $x\Rightarrow a$ iff
$f(x)\Rightarrow f^{-1}(a)$, or in usual terms, $f(x)\subseteq a$ iff
$f^{-1}(f(x))\subseteq f^{-1}(a)$.

The anti-diagonal universal $h_{x2}:f(x)\Rightarrow f^{-1}(f(x))$ is just the
truism $f^{-1}(f(x))\subseteq f^{-1}(f(x))$ or, equivalently, $f(x)\subseteq
f(f^{-1}(f(x)))$. The other universal $e_{a1}:f(f^{-1}(a))\Rightarrow
f^{-1}(a)$ is $f^{-1}(f(f^{-1}(a)))\subseteq f^{-1}(a)$ or, equivalently, the
other truism $f(f^{-1}(a))\subseteq f(f^{-1}(a))$. One over-and-back identity is:

\begin{center}
$f(x)\overset{h_{x2}}{\Rightarrow}f^{-1}(f(x))\overset{e_{f(x)}}%
{\Longrightarrow}f(x)$ iff $f(x)\subseteq f(f^{-1}(f(x)))\subseteq f(x)$ iff
$f(x)\subseteq f(x)$,
\end{center}

\noindent i.e., $f(f^{-1}(f(x)))=f(x)$. The other over-and-back identity is:

\begin{center}
$f^{-1}(a)\overset{h_{f^{-1}(a)}}{\Longrightarrow}f(f^{-1}(a))\overset{e_{a1}%
}{\Longrightarrow}f^{-1}(a)$ iff $f^{-1}(a)\subseteq f^{-1}(f(f^{-1}%
(a)))\subseteq f^{-1}(a)$ iff $f^{-1}(a)\subseteq f^{-1}(a)$,
\end{center}

\noindent i.e., $f^{-1}(f(f^{-1}(a)))=f^{-1}(a)$.

In the case of the limit and colimit adjunctions, the constant functor
$\Delta$ from $\mathbf{Set}$ to $\mathbf{Set}^{D}$ had both a right adjoint
($Lim$) and a left adjoint ($Colim$) so there was two-way determination
between sets and diagram functors. In the present case, the inverse image
functor $f^{-1}()$ has a left adjoint so the question arises of it having a
right adjoint. If it had a right adjoint, then that adjunction would be about
the determination from a subset $a$ in $\mathbf{A}$ to a subset $x$ in
$\mathbf{X}$ where the determinative relation $a\Rightarrow x$ holds if for
all $\alpha\in a$, $f^{-1}(\{\alpha\})\subseteq x$. A right adjoint applied to
an $\mathbf{X}$-object $x$ would give the $\mathbf{A}$-object that represents
all the possible determinees of $x$. The determinees of $x$ are its subsets of
the form $f^{-1}(a)$ for some $a$ in $\mathbf{A}$ so the representation of all
the possible determinees as an $\mathbf{A}$-object would be obtained as the
union $\cup\{a\mid f^{-1}(a)\subseteq x\}$ which might be denoted as $f_{\ast
}(x)$ and which can also be defined directly as: $f_{\ast}(x)=\{\alpha
\in\mathcal{A}\mid f^{-1}(\{\alpha\})\subseteq x\}$. This yields the
adjunction or Galois connection: $f^{-1}(a)\subseteq x$ iff $a\subseteq
f_{\ast}(x)$. The two Galois connections give two determinations through
universals in both directions between $\mathbf{X}$ and $\mathbf{A}$.

\section{Conclusions}

\subsection{Summary of Theory of Adjoints}

We have approached the question of \textquotedblleft What is category
theory?\textquotedblright\ by focusing on universal mapping properties and
adjoint functors which seem to capture much of what is important and universal
in mathematics. We conclude with a brief summary of the theory and with some
philosophical speculations.

There is perhaps some irony in the theory presented here to explain one of the
central concepts in category theory, the notion of an adjunction. We were
required to reach outside the conventional ontology and to acknowledge the
object-to-object chimera morphisms or heteromorphisms between categories. The
adjunction representation theorem, which shows that all adjunctions arise from
birepresentations of het-bifunctors of chimera morphisms was based on the
adjunctive square construction. The adjunctive square is a very convenient
and, indeed, natural diagram to represent the properties of an adjunction. An
adjunction couples two categories together so that there is a form of
determination that goes from the $x$-pair $\widehat{x}=(x,Fx)$ to an $a$-pair
$\widehat{a}=(Ga,a)$. In a commutative adjunctive square, the main diagonal
$(f,g)$ gives the pair of determinations $f=g^{\ast}:x\rightarrow Ga$ and
$g=f^{\ast}:Fx\rightarrow a $. All adjunctions can be obtained as the
birepresentation of the cross-category determinations expressed by a
het-bifunctor $Het:\mathbf{X}^{op}\times\mathbf{A}\rightarrow\mathbf{Set}$.
\ For many of the adjunctions encountered by working mathematicians, the
determinations can be specified concretely by some chimera morphisms or
heteromorphisms $x\Rightarrow a$. But the abstract het-bifunctor, where
$Het(\widehat{x},\widehat{a})=\{\widehat{x}=(x,Fx)\overset{(f,f^{\ast}%
)}{\longrightarrow}(Ga,a)=\widehat{a}\}$, reproduces (up to isomorphism) any
given adjunction as its birepresentation. In the adjunctive square format
(abstract or concrete), each determination factors uniquely through sending
and receiving universals at each end of the determination (the zig-zag factorization).

A powerful philosophical theme, which connects self-determination and
universality, emerged in the conceptual analysis of the role of certain
identity morphisms in each universal of an adjunction. Each universal is
self-participating but, in addition, a type of `self-determination' is
involved in the construction of each universal itself. The self-determination
was represented by the identity maps $1_{Fx}$ and $1_{Ga}$, and by their
associated universal maps in the isomorphisms of bifunctors:%
\[
Hom(Fx,a)\cong Het(x,a)\cong\mathcal{Z}(Fx,Ga)\cong Hom(x,Ga).
\]

\noindent In the zig-zag factorization, any determination given by an
adjunction can be factored indirectly using the anti-diagonal map to be
compatible with the self-determination represented by the universals on the
sending and receiving ends.

\subsection{Determination through universals: the main features}

In the adjunctions of category theory, we have seen that on the receiving end,
a determination might be expressed by a direct map $Fx\overset{g}{\rightarrow
}a$ or by an indirect map $Fx\overset{Fg^{\ast}}{\rightarrow}FGa$ factoring
through the receiving universal $FGa\overset{\varepsilon_{a}}{\rightarrow}a$.
In chimera terms, the same direct map $Fx\overset{g(c)}{\longrightarrow}a$ has
the over-and-back factorization through the indirect anti-diagonal map
$Fx\overset{z(c)}{\Longrightarrow}Ga$ and the chimera counit or receiving
universal $Ga\overset{e_{a}}{\Longrightarrow}a$. Mathematically the direct and
the indirect-through-the-universal determinations are equal but in the
empirical sciences there might be a question of whether a determinative
mechanism or process was of the first direct type or the second indirect type
factoring through the universal. With direct determination, the receiver has
the passive role of receiving the determination. In the second type of
mechanism, the receiver of the determination plays a more active or
self-determining role of generating a wide (`universal') range of
possibilities and then the determination takes place indirectly through the
selection of certain of those possibilities to be actively implemented.

Several main features of this determination through universals might be
singled out for the receiving case (the sending case is dual).

\begin{description}
\item[Universality:] While an external direct determination specifies or
determines a particular set of possibilities, the determination through a
universal constructs the object representing all the possibilities that might
be directly determined---as indicated by its universal mapping property.

\item[Autonomy:] The universal is constructed in a manner independent of any
external determiners (e.g., neither any $x$ nor any $f$ or $g$ were involved
in constructing the receiving universal in the conventional form
$FGa\overset{\varepsilon_{a}}{\rightarrow}a$ or in the chimera form
$Ga\overset{e_{a}}{\Longrightarrow}a$).

\item[Self-determination:] The morphism associated with the universal
potentially determines all the possibilities (e.g., the receiving universal in
either form was the adjoint correlate of the identity $1_{Ga}$).

\item[Indirectness:] The particularization comes only with the indirect factor
map that picks out or selects certain possibilities.

\item[Composite Effect:] The composition of the specific factor map followed
by the universal morphism then implements the possibilities to agree with the
given direct determination.
\end{description}

In grand philosophical terms, the factorization through universals of an
adjunction gives an approach, albeit in rather abstract mathematical terms, to
resolving what is perhaps the central conundrum of philosophy, the
reconciliation of external determination (\textquotedblleft
necessity\textquotedblright\ or \textquotedblleft heteronomy\textquotedblright%
) and self-determination (\textquotedblleft freedom\textquotedblright\ or
\textquotedblleft autonomy\textquotedblright). We turn to what this
determination through universals might mean in the life sciences and in social
philosophy.\footnote{In the physical sciences, there is the obvious
possibility of viewing the expansion of the wave packet in quantum mechanics
and then its collapse to realize a specific actuality as an application of
determination through universals. However, I am not prepared to investigate
that possibility here so the focus is on the life and human sciences.}

\subsection{Determination through universals in the life sciences}

The debate between selectionist and instructionist mechanisms can be seen in
this light. In the original debate about evolution, the Lamarckian position
was that the environment directly \textquotedblleft
instructed\textquotedblright\ the organism about adaptive features which were
then inherited by the organism's offspring. In the Darwinian selectionist
theory, the species generated a wide range of possibilities (e.g., through
mutations and sexual reproduction) and then the environment has only the
indirect role of selecting which features have a survival advantage and which
will thus be differentially propagated to the offspring.

The mathematics provides a highly abstract, atemporal, and idealized model so
one does not expect a perfect fit to any processes in the life or human
sciences. With that caveat, the main features of determination through
universals seem to be present in the selectionist account of biological evolution.

\begin{description}
\item[Universality:] The selectionist theory is an example of population
thinking because it is the population, not the individual organism, that
explores the universe of possibilities by variation through mutation and
sexual reproduction.\footnote{Some versions of Darwinian evolution have taken
the problem of the generation of variety more seriously than others. In
particular, Sewall Wright's shifting balance theory \cite{wright:evol} has
emphasized the advantages of having the population split up into various
subpopulations or `demes' that will encourage wider variation (a practice also
followed by artificial breeders). Separation of subpopulations allows more
variation to be tested and fixed but there also has to be migration between
the subpopulations so that `improvements' or `discoveries' will be transmitted
to the whole population. The proper mix is a question of shifting balances.}

\item[Autonomy:] In the modern treatment of genetic evolution, there is some
emphasis on the \textquotedblleft fundamental dogma\textquotedblright\ that
there is no information flow from the environment to `direct' the process of
generating genetic variety. This is the aspect of autonomy. The possibilities
are generated independent of any external determiners.

\item[Self-determination:] In the biological context, this is the
self-reproduction of organisms. In the mathematical example of the limit
adjunction, we saw how the determinees became determiners in the construction
of the limit $LimD$ so the determinees determined themselves via the
projection maps $e_{D}:LimD\Rightarrow D$. In a temporal context, the
switching of determinee and determiner roles would be sequential. This
switching of roles is involved in all biological reproduction, the offspring
becomes the parent. At the molecular level of reproduction, a DNA\ sequence is
first the determinee when it is formed on a given template, and then after the
splitting of the double helix, it becomes itself a template or determiner to
`determine' or reproduce itself.

\item[Indirectness:] The selective effect of the environment on the variety of
possibilities that have been generated is modelled mathematically by the
indirect factor map.

\item[Composite Effect:] The composite effect of this external selection and
potential self-reproduction of the variants is to `implement' or
differentially reproduce the selected variants. That is mathematically
modelled by the composition of the indirect factor map (which selects certain
possibilities) with the universal morphism---which `by itself' would
self-determine or reproduce all possibilities---so that the \textit{composite}
effect is to differentially amplify or `implement' the selected possibilities.
\end{description}

Today the idea of an instructive versus a selective process has been
generalized to a number of other processes. A common theme is that learning
processes are originally thought to be instructive (i.e., direct
$Fx\overset{g(c)}{\longrightarrow}a$) but are then found to be selective
(i.e., indirect through the universal $Fx\overset{Fg^{\ast}}{\rightarrow
}FGa\overset{\varepsilon_{a}}{\rightarrow}a$ or $Fx\overset{z(c)}%
{\Longrightarrow}Ga\overset{e_{a}}{\Rightarrow}a$). The key component of any
selectionist mechanism is the generator of diversity that generates the
`universal' range of possibilities (e.g., the construction in the mathematics
of going from $a$ to $FGa$ or simply to $Ga$) so that certain possibilities
can then be selected (e.g., by the indirect morphism $Fx\overset{Fg^{\ast}%
}{\rightarrow}FGa$ or $Fx\overset{z(c)}{\Longrightarrow}Ga$) to determine the
eventual outcome by a self-determinative process (e.g., by the canonical
morphism $FGa\overset{\varepsilon_{a}}{\rightarrow}a$ or $Ga\overset{e_{a}%
}{\Longrightarrow}a$).

A particularly striking application of the selectionist approach is to the
immune system. The early theories were instructional; the external molecule or
antigen would enter the system and instruct the immune mechanism with its
template to construct antibodies that will neutralize the antigen. In 1955,
Niels Jerne \cite{jerne:nst} proposed the selectionist theory of the immune
system which is now accepted (with variations added by many other
researchers). The immune system takes on the active role of generating a huge
variety of antibodies and the external antigen has the passive role of simply
selecting which antibody fits it like a key in a lock. Then that antibody is
differentially amplified in the sense of being cloned into many copies to
lock-up the other instances of the antigen. A similar example was the
originally instructivist account of bacteria `learning' to tolerate
antibiotics or to consume a new substance but now these processes are
recognized as selectionist. A wide variety of bacterial mutations are
constantly being generated and those that can tolerate antibiotics or digest a
new substrate will differentially thrive in such an environment.

Peter Medawar (who in addition to Jerne received a Nobel Prize for work
related to the immune system) illustrated the difference between an
instructive and selective (or \textquotedblleft elective\textquotedblright)
mechanism using as an analogy the difference between a phonograph (or
\textquotedblleft gramophone\textquotedblright) and a jukebox. With a
phonograph, one has to externally supply the specific musical instructions (a
record) to the machine which then plays the record. \ But a jukebox has a wide
repertoire of musical instructions (records) inside it; externally there is
only the selection of the record to be played.

The analogy with the determination through universals can be illustrated using
the conventional exponential adjunction. Taking the set $X=1$, a single record
played on the phonograph might be compared to the function $g:1\times
A\rightarrow Y$ where $A$ is a set of parameters---going over the set of
parameters $A$ `plays' the record to produce the music $Y$. But the jukebox
contains within it a large `universal' repertoire of records $Y^{A}$ which
might be played. The adjoint transpose $g^{\ast}:1\rightarrow Y^{A}$ picks out
the same record that played as $g:1\times A\rightarrow Y$ out of the universal
repertoire $Y^{A}$. The operation of the jukebox playing a record from its
internal repertoire is represented by the counit $\varepsilon_{Y}:Y^{A}\times
A\rightarrow Y$. To get the same effect with the jukebox as with the
phonograph, the functor $-\times A$ is applied to the selection $g^{\ast
}:1\rightarrow Y^{A}$ of the record to be played:

\begin{center}
$%
\begin{array}
[c]{ccc}%
"Phonograph" & 1\times A\overset{g}{\longrightarrow}Y=1\times A\overset
{g^{\ast}\times A}{\longrightarrow}Y^{A}\times A\overset{\varepsilon_{Y}%
}{\longrightarrow}Y & "Jukebox"
\end{array}
$\ \qquad\qquad\ \ \ \ \ \ \ 
\end{center}

\begin{quote}
\noindent{\footnotesize During the past ten years [1950s], biologists have
come to realize that, by and large, organisms are very much more like
juke-boxes than gramophones. Most of the reactions of organisms which we were
formerly content to regard as instructive are in fact elective.\cite[p.
90]{med:fut}}
\end{quote}

\noindent In the mathematics, the morphisms are equal, but empirically the
point is that a selective mechanism is represented by the factorization
through the universal, not by the direct morphism that mathematically would
give the same end results (e.g., the playing of the record).

In a selectional mechanism, there is also instruction or determination but it
comes from within---the autonomy and self-determination that emerges in the
mathematics with the counit $\varepsilon_{Y}:Y^{A}\times A\rightarrow Y$ being
the adjoint transpose of the identity map $1_{Y^{A}}:Y^{A}\rightarrow Y^{A}$
(the external $X$ plays no role). The point about the jukebox was not that it
had no musical instructions (records) but that they were embodied within the
entity rather than externally supplied. \textquotedblleft The instructions an
organism contains are not musical instructions inscribed in the grooves of a
gramophone record, but \textit{genetical} instructions embodied in chromosomes
and nucleic acids.\textquotedblright\ \cite[p. 90]{med:fut} Heinz Pagels makes
the same point connecting evolution and the immune system.

\begin{quote}
{\footnotesize Like evolution, the immune response is also a selective system
in which the system is instructed from within---the genetic instructions plus
random variation---but the selection depends on the external environment---the
specific invading antigens. \cite[p. 265]{pag:dream}}
\end{quote}

Language learning by a child is another example of a process that was
originally thought to be instructive. \ But Noam Chomsky's theory of
generative grammar postulated an innate universal grammar (the instructions
from within) that would unfold according to the linguistic experience of the
child. \ The child did not `learn' the rules of grammar; the linguistic
experience of the child would select how the universal mechanism would develop
or unfold to implement one rule rather than another. This connection is not
new. Niels Jerne's Nobel Lecture was entitled \textit{The Generative Grammar
of the Immune System}.

An everyday example of indirect determination is a person's understanding of
spoken language. \ The naive viewpoint is that somehow the meaning of the
spoken sentences is transmitted from the speaker to the listener. But, in
fact, it is only the physical sounds that are transmitted. The syntactic
analysis and the semantic component have to be generated by the listener so
the heard sounds only have the role of selecting which generative processes
will be triggered. \ Here again, Chomsky has emphasized the universality of
the internal mechanism to generate an understanding of a potential infinity of
sentences which have never been heard before. Descartes emphasized this
universality of language and reason: \textquotedblleft reason is a universal
instrument which can serve for all contingencies\textquotedblright\ \cite[p.
116]{des:pw1} so Chomsky has referred to the generative grammar approach as
\textquotedblleft Cartesian linguistics.\textquotedblright\ 

\begin{quote}
{\footnotesize In summary, one fundamental contribution of what we have been
calling "Cartesian linguistics" is the observation that human language, in its
normal use, is free from the control of independently identifiable external
stimuli or internal states and is not restricted to any practical
communicative function, in contrast, for example, to the pseudo language of
animals. It is thus free to serve as an instrument of free thought and
self-expression. The limitless possibilities of thought and imagination are
reflected in the creative aspect of language use. The language provides finite
means but infinite possibilities of expression constrained only by rules of
concept formation and sentence formation, these being in part particular and
idiosyncratic but in part universal, a common human endowment. \cite[p.
29]{chom:cl}}
\end{quote}

\noindent The general features of universality, autonomy (independence from
external stimulus control), and self-determination are clear.

After Gerald Edelman received the Nobel prize for his work on the selectionist
approach to the immune system, he switched to neurophysiology and developed
the theory of neuronal group selection or neural Darwinism.

\begin{quote}
{\footnotesize [T]he theoretical principle I shall elaborate here is that the
origin of categories in higher brain function is somatic selection among huge
numbers of variants of neural circuits contained in networks created
epigenetically in each individual during its development; this selection
results in differential amplification of populations of synapses in the
selected variants. In other words, I shall take the view that the brain is a
selective system more akin in its workings to evolution than to computation or
information processing.\cite[p. 25]{edel:nd}}
\end{quote}

\noindent In \noindent simpler terms, the brain generates an immense variety
of groups of neural circuits (like the variety of antibodies) and then
external stimuli only select which neural groups will be differentially
amplified. What at first looks like the external environment instructing the
brain is seen instead as a selectionist process.

From antiquity, some schools of thought (e.g., Neo-Platonism) have emphasized
the general point that understanding and learning are mistakenly seen as a
direct instructive process rather than an indirect composite effect of
catalyzing a universal internal generative process. In the early fifth
century, Augustine in \textit{De Magistro }(The Teacher) made the point.

\begin{quote}
{\footnotesize But men are mistaken, so that they call those teachers who are
not, merely because for the most part there is no delay between the time of
speaking and the time of cognition. And since after the speaker has reminded
them, the pupils quickly learn within, they think that they have been taught
outwardly by him who prompts them.(Chapter XIV)}
\end{quote}

\noindent Wilheim von Humboldt made the same point even recognizing the
symmetry between speaker and listener.

\begin{quote}
{\footnotesize Nothing can be present in the mind (Seele) that has not
originated from one's own activity. Moreover understanding and speaking are
but different effects of the selfsame power of speech. Speaking is never
comparable to the transmission of mere matter (Stoff). In the person
comprehending as well as in the speaker, the subject matter must be developed
by the individual's own innate power. What the listener receives is merely the
harmonious vocal stimulus. \cite[p. 102]{hum:herm}}
\end{quote}

\noindent In the mathematics of adjunctions we have seen the symmetry between
the self-determination involved in both the sending and receiving universals
(e.g., the zig-zag factorization).

In all these cases---from the beginnings of evolutionary selection up through
the highest human functions---one of the central points is that being on the
receiving end of a determination does not imply passivity. There is an
indirect mode of determination through a receiving universal where the
receiver in some form actively generates or has the capacity to actively
generate a well-nigh `universal' set of possible receptions (received
determinations or determinees). \ Then the determination (e.g., learning)
through the receiving universal takes the form of selecting which `message' is
actively generated so that receiving the determination is quite consistent
with the active self-determination of the receiver. Such a mechanism seems key
to understanding how an organism can perceive and learn from its environment
without being under the direct stimulus control of the environment---thus
resolving the ancient conundrum of receiving an external determination while
exercising self-determination.

In the adjunctions of category theory, we have seen an abstract version of
this conceptual structure of indirect determination through universals that
express self-determination. The importance of these conceptual structures
\textit{in mathematics} was the original reason why we have focused on
adjunctions. We have given a few hints at how these adjunctive structures of
determination through universals might also model processes of central
importance in the life sciences that serve, in varying degrees, to make
external determination more indirect and thus more consistent with self-determination.

\subsection{Determination through universals in social philosophy}

The conceptual structure of determination through universals might also be
used as a normative model. Here the \textit{locus classicus} is Immanuel Kant.
The mathematics of adjunctions sustains Kant's philosophical intuition that
universality (the first version of the categorical imperative) is closely
related to autonomy and self-determination (the second and third versions of
the categorical imperative). However, it is not clear that he worked out the
correct notion of universality that would correspond to, say, the autonomy
principle of always treating persons as ends-in-themselves rather than just as
means (\cite{kant:ground}; see Chapter 4 in $\cite{ell:it}$). In any case,
there may be something to Michael Arbib and Ernest Manes' linking of
\textquotedblleft The Categorical Imperative\textquotedblright\ (subtitle of
their book) to universal constructions in category theory \cite[p.
vii]{arbib:asf}. And it might be noted that Kant also saw an active role for
the mind in perception and cognition somewhat along the lines described in the
previous section.

Perhaps the most important applications come in political and economic theory.
In political theory, all heteronomous governance relations can be restructured
to form a political democracy so that people are, at least in theory, jointly
self-governing. In economics, the problematic determinative relation is that
of the employer and employee. But it is also always possible to restructure a
firm as a workplace democracy so that everyone working in it is jointly
`self-employed' or self-managing.{\footnote{The relevant political theory is
outlined in \cite{ell:tvc} and the concepts of workplace democracy are dealt
with at book length in $\cite{ell:pc}$ (which can be downloaded from
www.ellerman.org ).}}

More generally, normative questions of heteronomy versus autonomy arise in
social philosophy where determinative relations between persons are the
central topic. The setting is where one person (or group of persons) has the
sender or `determiner' role and another person (or group of persons) is in the
receiver or `determinee' role. For the person in the sending role trying to
influence, instruct, counsel, control, or determine others, the message is
that this can structured as firstly being self-determination---leading by
example, practicing what one preaches, and teaching something by doing it
oneself. But the greater problem is the self-determination or autonomy of the
person or persons in the receiving position. The mathematical analogy shows
that (in an adjunctive context) there is always a way to rearrange matters so
that any external determination becomes indirect by factoring through the
receiving universal that realizes the self-determination of the receiver.

Across human affairs, there are relationships of teacher to student, manager
to subordinate, counselor to client, psychologist to patient, and helper to
doer where the first party tries to influence, control, or otherwise determine
the actions and beliefs of the second party. The perennial conundrum is that
most `help' or educational instruction occurs in a manner that overrides or
undercuts the self-determination and autonomy of the persons in the receiving
position. It is a most subtle matter to see how such heteronomous external
determination might be `factored' to take an indirect form that would respect
the self-determination and autonomy of the people in the receiver
role.{\footnote{These matters have been dealt with at book length in
$\cite{ell:hpht}$ (with some emphasis on development assistance).}}

This theory of adjoints gives an account of an adjunction based on
determination through universals expressing a type of self-determination at
both the sending and receiving end of a determination. The salient features
were summarized in adjunctive squares which show that---by the zig-zag
factorization---it is always possible to factor any determination expressed by
an adjunction indirectly through the sending and receiving universals. It is
heartening to see the \textquotedblleft unreasonable effectiveness of
mathematics\textquotedblright\ (to use Eugene Wigner's expression) to capture
this basic conceptual theme of structuring an external determinative
relationship to be indirect and compatible with self-determination on each end
of the determination.

\section{Acknowledgment}

I am indebted to Steve Awodey, Vaughan Pratt, Colin McLarty, and John Baez for
comments and suggestions about the many earlier versions of this paper.

\medskip

\noindent Department of Philosophy

\noindent University of California at Riverside

\noindent e-mail: david@ellerman.org

\noindent web-site: www.ellerman.org

\end{document}